\documentclass[12pt]{article} 
\usepackage{verbatim,amsmath,amssymb}

\newcommand{\rf}[1]{(\ref{#1})}

\newcommand{\numer}[1]{\label{#1}}

\newcommand{\cd}{{\cal D}}
\newcommand{\cc}{{\cal C}}
\newcommand{\cb}{{\cal B}}
\newcommand{\cu}{{\cal U}}
\newcommand{\ch}{{\cal H}}
\newcommand{\cf}{{\cal F}}

\usepackage{epsfig}
\usepackage{graphicx}
\usepackage{index}
\usepackage{amscd}
\input math111.def
\input buch.def
\numberwithin{equation}{section}
\newtheorem{conjecture}[theor]{Conjecture}

\def\span{\mbox{\rm span }}

\newif\ifmarglab
\makeatletter
\def
\label#1{\@bsphack\ifmarglab\marginpar{LAB:#1}\fi\if@filesw {\let\thepage\relax
   \def\protect{\noexpand\noexpand\noexpand}%
   \edef\@tempa{\write\@auxout{\string
      \newlabel{#1}{{\@currentlabel}{\thepage}}}}%
   \expandafter}\@tempa
   \if@nobreak \ifvmode\nobreak\fi\fi\fi\@esphack}
\makeatother

\begin{document}
\title{Rearrangements with supporting Trees, Isomorphisms\\
and Combinatorics of coloured dyadic Intervals
\footnote{Partially supported by SPADE2 programm at IM PAN. Paul F.X. M\"uller supported by FWF P 20166-N18.} 
}
\author{Anna Kamont, Paul F. X . M\"uller}
\maketitle

\begin{abstract}
We determine a class of rearrangement operators acting on dyadic intervals that admit a supporting tree. 
This condition implies that the associated rearrangement operator has a bounded vector valued extension to $L^p_E$, 
where $E$ is a UMD space.

We prove the existence of a large subspace $X_p\subset L^p$ on which a bounded rearrangement operator acts as an isomorphism.

Moreover, we study winning strategies for a combinatorial two person game played with coloured collections of dyadic intervals.
\end{abstract}



\tableofcontents
\newpage
\section{Introduction}
We study the connections  between rearrangement operators of the Haar system and martingale transforms. We restrict to rearrangements $\t $ acting on 
dyadic intervals such that $|\t(I)| = |I|  $ and operators 
given by 
$$ T(h_I)= h_{\t(I)} . $$
Our  chief interest lies in obtaining workable and directly verifyable 
descriptions of those rearrangements $\t$ for which $T \otimes \Id_E $
is bounded on $L^p_E$  for all Banach spaces  $E$  in the UMD class.

We isolate first a purely combinatorial condition on $\tau$  that implies the 
boundedness of $T \otimes \Id_E $  on $L^p_E $ for $1 < p < \infty,$
and $E$ in the UMD class.Thereby we convert    the analytic question into a
combinatorial problem. We let $\cF$ be a collection of dyadic intervals 
and $\{ A_I :\, I \in \cF \} $ be a tree of measurable sets such that 
$$ A_I \cap A_J \ne \es \quad\text{implies}\quad A_I \sbe A_J \quad\text{or}\quad A_J \sbe A_I .$$ 
We say that $\t : \cF \to \t(\cF) $ admitts the supporting tree 
$\{ A_I :\, I \in \cF \} $ if there exists $c > 0 ,\delta >0$ so that $ |A_I | \le c |I| $
and  
\begin{equation}\label{17jan91}
   | I \cap A_I| \ge  \d |I|   \quad\text{and}\quad  | \t(I) \cap A_I| \ge  \d |I| , 
\quad \quad I \in \cF.
\end{equation}
If   $\t : \cF \to \t(\cF) $ admitts a supporting   tree, 
then, by an application of Stein's martingale inequality, 
$T \otimes \Id_E $ is bounded on 
$$ \span _{L^p_E} \left\{h_I x_I : \, x_I \in E ,\, I \in \cF \right\} ,\quad
1 < p < \infty . $$
Hence if we ask for the $L^p_E $ boundedness of $T \otimes \Id_E $ the ensueing combinatorial problem consists 
in  decomposing the collection of all dyadic intervals into 
$ \cF_1 ,\dots ,\cF_N $ so that the restrictions 
$\t : \cF_i \to \t(\cF_i ) $ admit a supporting   tree.
For such a decomposition to exist it is 
necessary  that 
$T$ (scalar valued) 
is an isomorphism
on some $L^p  $ with $1 <p \ne 2 < \infty . $  The results of this paper are
all related to the open problem  
whether this  
is also a {\em sufficient condition.}
That is, we are concerned with the following extension problem for 
rearrangement operators:
\begin{conjecture}\numer{con1}
Assume that  the scalar valued rearrangement operator
$ T $ is an isomorphism on $L^p, \, p \ne 2 .$
Is it true that $T \otimes \Id_E $ is an isomorphism on $L^p_E , $
for  $ 1 < p < \infty$ and any $E$ in the UMD class.
\end{conjecture}
The vector valued extension problem for rearrangement operators
as formulated above is meaningful only within the class of
{\it isomorphisms} on $L^p.$ Indeed \cite{gmpre} and \cite{pfxmii}
contain examples of {\it bounded} rearrangement operators on
$L^p , 1 <p <2 $  and asssociated examples of UMD spaces $E$ 
so that the vector valued extension
$T \otimes \Id_E $ is {\it not} a   bounded operator on $L^p_E . $

Semenov's theorem \cite{ns} 
provides an intrinsic criterion   for  $ \t$ 
so that $T$  is an $L^p$ isomorphism.  
Thus by considering $T$  {\em together with its }  
vector valued extension 
$ T\otimes \Id_E $ we are led to the following combinatorial problem concerning 
tree structures and rearrangements. 
\begin{conjecture}\numer{con2}
Does the hypothesis 
$$ \left| \bigcup_{J\sbe I } \t(J) \right| \le C|I| \quad\text{and}
\quad \left| \bigcup_{J\sbe I } \t^{-1} (J) \right| \le C|I|,$$ 
 for any dyadic $I$ imply, that the entire collection of dyadic intervals 
can be decomposed into $\cF _ 1, \dots , \cF_N ,$ with $N=N(C),$  
so that the restrictions 
$$\t : \cF_i \to \t(\cF_i ) \quad\text{admit a supporting tree?}$$ 
\end{conjecture}
In this paper we give  partial solutions to the conjectures~\ref{con1} and \ref{con2}.
Our results are connected to a wider set of problems as follows:
\begin{enumerate}
\item Seeking understanding of scalar valued operators by studying simultaneously their vector valued extensions is a central line of 
investigation in Banach space theory. Classical and authoritative accounts thereof are
 \cite{tficm} by  T. Figiel and  \cite{pelczicm} by A. Pe{\l}czy{\'n}ski.
In the context of rearrangement operators, the search for supporting trees is just motivated by
our attempts to prove the boundedness of vector valued rearrangements.

\item
Section 4.7. in \cite{pisier-2008} by G. Pisier contains the question of 
describing
the class of $L^p$ bounded operators $R$ so that 
 $R\otimes \Id_E$
is bounded on $L^p_E ,  $  for any UMD space $E.$ 
 Calderon-Zygmund integrals and martingale transforms share this property.
The theorems in the present paper aim at a description of the 
rearrangement operators in this class.
\item The rearrangement problems  arising  with the unconditionality of the Franklin system and its
generalisations \cite{figbases} , \cite{km} are special cases respectively 
model cases for the problems treated  in Section~\ref{sem.s1}.
\item The recent extrapolation theorems \cite{gmpre}, \cite{pfxmii}, for  bounded and invertible  vector valued 
rearrangement operators on the  $L_E^p \, (1 <p < \infty) $ scale provide
the structural support for the conjectures~\ref{con1} and \ref{con2}. 
Indeed, if the latter were true 
the former extrapolation results would simply follow from the well known scalar case.
\item The examples in \cite{gmpre} and \cite{pfxmii}
of $L^p$ bounded rearrangement operators $T$ 
for which
$$  \left\|  T\otimes \Id _E \right\| _{L_E^p} = \infty \quad\text{and}\quad 
\left\|  T^{-1}\otimes \Id _E \right\|_{L_E^p} =\infty,$$
provide the motivation for proving the subspace theorem in Section~\ref{sub}.
Depending on $\tau$ we determine a large subspace $X_p\sbe L^p$
on which nevertheless  the restricted operator $T_{|_{X_p}}\otimes \Id _E $
acts as an isomorphism so that,
 $$  \left\|  T_{|_{X_p}}\otimes \Id _E \right\| _{L_E^p}\cdot 
\left\|  T^{-1}_{|_{T(X_p)}}\otimes \Id _E \right\|_{L_E^p} <\infty, 
$$
for $1 < p < \infty$ and any UMD space $E .$
\end{enumerate}
The organisation of the paper is as follows. In Section
\ref{prel} 
we review the concepts and theorems used in the paper. 
In Section 
\ref{sub}
we prove a subspace theorem for rearrangement operators.
For an $L^p$ bounded rearrangment operator $T$ we determine a block basis 
$\{\widetilde{h_I}\} ,$ equivalent to the entire Haar system,
such that $T$ acts as an isomorphism on the subspaces spanned by 
$\{\widetilde{h_I}\} .$ This done by constructing a tree that simultaneously
supports $\{\widetilde{h_I}\} $ and its image $\{T(\widetilde{h_I})\}. $
Section \ref{sem.s1} treats  special shift operators $S. $  
We show that $S\otimes \Id_E $ is $L^p_E$ bounded for $E$ in the UMD class,
provided that the associated shift parameters form a decomposable sequence
in the sense of Defintion~\ref{sem.d1}. In Section~\ref{c.s1}
we study thoroughly the problem of finding $(d, \eta)$ homogeneous 
decompositions of a collection $\cC$ of dyadic intervals, 
that preserves a pre-existing $(d,\eta)$ homogenous decomposition
of a fixed subcollection $ \cC' $ in $ \cC .$ We obtain conditions 
for its existence and determine examples for which this problem is 
without solution. While -we think- this is a combinatorial problem 
interesting in itself we present it here since it should 
support the construction of trees for rearrangements. 

\section{\numer{prel} Preliminaries}
\paragraph{Dyadic intervals and the Haar system.}
See e.g. \cite{m}.
We let $\cD $ denote the collection of (half-open) dyadic intervals
contained in [0,1]. Thus
$$
\cD=\{
[(k-1)2^{-n} , k 2^{-n}[ , \quad 1\le k \le 2^{n} , \quad n \in \bN \}. $$
For $n \in \bN $ write  
$\cD_n = \{ I \in \cD : |I| = 2^{-n} \}. $  
For a collection of dyadic intervals $\cE $ we use the $*-$ notation to 
denote the pointset covered by $\cE $ thus
$$ \cE^* = \bigcup_{I \in \cE } I . $$
For $I \in \cD $ denote by $Q(I)$ the collection of all dyadic intervals that 
are contained in $I . $ 

Denote by $ \{ h_I :  I \in \cD \} $
the $L^\infty-$ normalized Haar system, where $h_I$ is supported on $I$ 
and 
$$
h_I = \begin{cases} 1 \quad &\text{ on the left half of }I ;\\ 
                    -1  \quad &\text{ on the right half of }I . 
      \end{cases}
$$ 
The Haar system is an unconditional basis in $L^p,$ $(1 < p < \infty) .$
For $ f \in  L^p$ we define its dyadic square function as 
$$ S(f) = \left( \sum_{I \in \cD } \la f,\frac{h_I}{|I|}\ra ^2 1_I \right)^{1/2}.$$
The Marcinkiewicz \cite{jm37} interpretation of  R.E.A.C. Paley's theorem \cite{rp32}
asserts that 
$$ c_p \| f\|_{L^p} \le \| S(f)\|_{L^p} \le C_p \| f\|_{L^p} , $$
where $C_p \sim p^2/ (p-1)$ and $c_p  = C_p^{-1} . $
Recall also the definition of dyadic $H^1,$ defined by the relation 
$ f \in H^1$ if   $S(f) \in L^1$   and
  $$  \| f\|_{H^1} = \| S(f) \|_{L^1} . $$
\paragraph{The Haar system in Bochner-Lebesgue spaces.}
\cite{fw}, \cite{ma2}.
Let $ 1 < p <\infty .$
For a Banach space $E$ we denote by $L_E^p ,$
the Bochner Lebesgue space of E valued and $p -$ integrable functions
on the unit interval. 
We say that a Banach space $E$ satisfies the UMD property if 
there exists $C_p > 0 $  so that for any finite collection 
$x_I \in E $ 
$$
\left \| \sum \pm x_I h_I \right\|_{ L_E^p}
\le C_p \left \| \sum  x_I h_I \right\|_{ L_E^p} . $$
\paragraph{Kahane's principle of contraction and Kahane's inequality.}
\cite{kahane},
\cite{w}.
Let $\{r_n \} $ denote the sequence of independent 
$\{+1, -1\}$ valued Rademacher functions. Let 
$x_n \in E $ be a sequence in a Banach space $E$
and let $ a_n \in \bR $ so that $ |a_n| \le 1 . $
Then, , 
$$\int_0^1 \| \sum _{n=1}^N  r_n(t) a_n x_n \|_E dt \le 
\int_0^1 \| \sum _{n=1}^N  r_n(t)  x_n \|_E dt . 
$$ 
 where $C >0 $ is independent of $N.$
We apply the above principle of contraction in combination with
the  Kahane's inequality asserting that 
$$
\left(\int_0^1 \| \sum _{n=1}^N  r_n(t)  x_n \|_E ^pdt \right)^{1/p}\le 
C_p
\int_0^1 \| \sum _{n=1}^N  r_n(t)  x_n \|_E dt, \quad  1<p<\infty . 
$$ 
\paragraph{A Martingale Inequality.} 
\cite{m} \cite{fw}.
The following vector valued martingale inequality is due to 
Bourgain and known as 
Bourgain's version of E.M. Stein's martingale inequality. 
It asserts the following. Let $E$ be a Banach space with the 
UMD  property. For any increasing sequence of $\s-$algebras
$\cF_n$ and any sequence  $ f_n \in L^p_E $ with $1 < p <\infty , $
$$
\int_0^1 \| \sum_{n=1}^N r_n(t) \bE (f_n | \cF_n ) \|_{ L^p_E}^p dt 
\le C_p \int_0^1 \| \sum_{n=1}^N r_n(t)  f_n  \|_{ L^p_E}^p dt. 
$$
where   $\bE (f | \cF )$ denotes conditional expectation with respect to 
$\cF$ and where  as above $\{r_n\}$
 denotes the sequence of independent Rademacher functions. 
\paragraph{The Theorem of Mazur.}
\cite{w}. Let $E $ be  Banach space and $ ( x_n ) $ be a sequence
in $E$ with weak limit $x\in E . $ Then there exist
a sequence $a_{n,j} \in [0,1]$ $ N_j \in \bN $ so that 
$$ \sum_{n = 1}^{N_j}  a_{n,j} = 1  \quad 
\text{ and } \quad 
 y_j =  \sum_{n = 1}^{N_j}  a_{n,j}x_n , \quad \quad j \in \bN , $$
converges to $x$ in norm of $E,$ that is, 
$ \lim_{j\to \infty} \|y_j -x \|_E = 0 . $

\paragraph{Semenov's Theorem. } 
\cite{ns}.
 Let  $\t : \cD \to \cD $ be  bijective satisfying
\begin{equation}\numer{10okt1}
 |\t(I)| = |I| , \quad I \in \cD .
\end{equation} 
The induced 
rearrangement operator is the linear extension of the map
$$ T : h_I \to h_{\t(I)}. $$ 
Let  $1 < p < 2 . $ Semenov's theorem asserts that 
$T$ is bounded on $L^p  $ if there exists $C>0$ so that  
 \begin{equation}\numer{10okt2} 
 |\t (Q(I))^*|\le C|I|,  \quad I \in \cD . 
\end{equation}
 Moreover, condition \eqref{10okt2} characterizes the boundedness of 
$T$ on $H^1. $
Specializing further we consider shifts defined 
by 
$$ \tau (I) = I + |I|m(|I|)  ,\quad I  \in\cD . $$
For shifts, Semenov's condition specializes as follows:
Shift operators are bounded on $L^p$ if there exists $K$ 
so that the set $ \t (Q(I))^* $ can be covered by $K$ 
dyadic intervals of the same length as $I.$
The best known of these shifts are those used by 
T. Figiel \cite{figbases, figsingular} to  study of Calderon-Zygmund 
operators,
$$ \tau_m (I) = I + m|I| ,\quad I  \in\cD . $$
Their norm estimates are given by T.Figiel's theorem \cite{figbases, figsingular}.
Below we apply it for fixed and  small values of $m .$  
\begin{theor}
\numer{tf.t1}
 The linear extension of   $T h_I = h_{\tau_m(I)}$
defines  a bounded operator on $L^p_E\,(1<p<\infty)$ for each Banach space $E$ with UMD property,
and 
$$
\| T: L^p_E \to L^p_E \| \leq C_p(E) (1 + \log m).
$$
\end{theor}

\paragraph{Dyadic trees.} Let $\cD$ be a collection of dyadic intervals
We say that  $\{E_I: I \in \cD\}$  is  a dyadic tree of sets
if, the following conditions hold: 
\begin{enumerate}
\item   Each of the sets $E_I$ is a finite union of dyadic intervals. 
\item There exists $C>0 $ so that 
\begin{equation}  |I|/C \le | E_I| \le  C |I| .
\end{equation} 
\item If $I_1$  is the left half of  $I $ and 
 $I_2$ is its  right  half 
  then 
\begin{equation}  E_{I_1} \cup   E_{I_2} \sbe   E_{I}  \quad\quad\text{and}   \quad\quad  E_{I_1} \cap   E_{I_2}
= \es . \end{equation} 
\end{enumerate}  
 Let $\{E_I: I \in \cD \}$  be  a dyadic tree and $\a \in \bR .$ Then we also 
use the term dyadic tree for the translates $\{\a + E_I: I \in \cD \}.$
We say that a collection of measurable sets $\{H_I: I \in \cD \}$ 
supports a dyadic tree  
$\{E_I: I \in \cD \}$
if there exists $ \d > 0 $ so that  for $ I \in \cD $
\begin{equation} 
H_I \sbe E_I \quad\quad \text{and}  \quad\quad|H_I| \ge \d |E_I|.
\end{equation} 
\paragraph{Trees and nested collections.}
Let $\cF $ be a subset of all  dyadic intervals.
 We say that $\{A_I: I \in \cF \}$
 is a tree  (or equivalently a nested collection)
of measurable sets if for $I, J  \in \cF , $
\begin{equation}  
A_I \cap A_J \ne \es  \quad\text{implies that}\quad 
 A_I \sbe A_J \quad \text{or} \quad A_J \sbe A_I . 
\end{equation} 
\paragraph{Rearrangements with supporting trees.}
Recall that the rearrangements we consider satisfy $|\t(I)| = |I| . $
We say that  $\t : \cF \to \cD $ admitts the supporting 
tree  $\{A_I: I \in \cF \}$
if there exists $C > 0 ,\delta >0$ so that 
$$  |A_I | \le C |I| ,$$ 
 
$$  | I \cap A_I| \ge  \d |I|   \quad\text{and}\quad  | \t(I) \cap A_I| \ge  \d |I| , 
\quad \quad I \in \cF.$$
The interest in the notion of rearrangements  admitting a supporting tree
comes from the following observation, obtained by merging 
\cite[Proposition 6.9]{jmst} with \cite{figbases, figsingular}.
\begin{theor}
\numer{sem.p20}
Let $\tau: \cF\to \cD$ be a rearrangement 
admitting  a supporting tree (with constants $C>0, \d > 0.$)
Let 
$$ X_p = {\rm span}_{L^p_E}\{x_I h_I, I \in \cF, x_I \in E \} \quad\text{and} \quad
Y_p = {\rm span}_{L^p_E}\{x_I h_{\tau(I)}, I \in \cF, x_I \in E \}.$$
Then $$T_\tau h_I = h_{\tau(I)}$$ extends to an isomorphism 
$$  T_{|X_p} \otimes \Id _E: X_p \to Y_p $$
so that 
$$
\left\| T_{|X_p} \otimes \Id _E\right\| \cdot \left\|  T^{-1}_{|Y_p} \otimes \Id _E\right\| \le f(C , \d, E ) .
$$
\end{theor}
\proof The proof is based on the contraction principle and 
Stein's martingale inequality.  
By hypothesis  $\t : \cF \to \cD $ has a supporting tree, 
say $\{ A_I : I \in \cF \} .$  Let $N \in \bN  $ and 
define two families of 
increasing $\s -$algebras,
 $$
\cA_N = \s \{ A_I :   |I| = 2^{-N} \}  
\quad\text{and}\quad
\cF_N = \s \{ I :   |I| = 2^{-N} \}. 
$$
We translate the hypothesis into pointwise estimates 
for conditional expectations. Let $I \in \cF $ with $|I| = 2^{-N } ,$ then
$$ 1_{\t(I) } \le \d^{-1}\bE( 1_{A_I \cap \t (I) }| \cF_N ) 
\quad\text{and}\quad 
1_{A_I \cap \t (I) } \le 1_{A_I} . $$
Now fix $x_I \in E $  for $I \in \cF.$  For $I \in \cD \sm \cF $ put $x_I = 0 .$ With the UMD property on 
$E ,$ the contraction principle,
and Bourgain's version of Stein's Martingale inequality we get, 
\begin{equation}
\begin{aligned}
\int_0^1 
\| \sum_{N \in \bN } \sum_{I \in \cD_N} r_I(t) x_I  1_{\t(I) }\|_{L^p_E}dt
&\le 
 C_p \d^{-1} \int_0^1\| \sum_{N \in \bN } \sum_{I \in \cD_N} r_I(t) x_I  \bE( 1_{A_I \cap \t (I) }| \cF_N )\|_{L^p_E}dt\\
&\le 
 C_p \d^{-1} \int_0^1\| \sum_{N \in \bN } \sum_{I \in \cD_N} r_I(t) x_I  1_{A_I  }
\|_{L^p_E}dt.
\end{aligned}
\end{equation}
Next exploit the second part of the hypothesis
$$ 1_{A_I } \le C\d^{-1}\bE( 1_{A_I \cap I }| \cA_N )  
\quad\text{and}\quad 
1_{A_I \cap I } \le 1_{I} , $$
and continue  again with  the contraction principle and the 
martingale inequality.
\begin{equation}
\begin{aligned}
\int_0^1\| \sum_{N \in \bN } \sum_{I \in \cD_N} r_I(t) x_I  1_{A_I  }
\|_{L^p_E}dt
&\le 
 C_p \d^{-1} \int_0^1\| \sum_{N \in \bN } \sum_{I \in \cD_N} r_I(t) x_I 
\bE( 1_{A_I \cap I }| \cA_N )
\|_{L^p_E}dt\\
&\le 
 C_p \d^{-1} \int_0^1\| \sum_{N \in \bN } \sum_{I \in \cD_N} r_I(t) x_I 
  1_{ I} 
\|_{L^p_E}dt.
\end{aligned}
\end{equation}
The reverse estimate follows in a similar fashion.
\endproof

\section{\numer{sub} Subspace Theorems for Rearrangement Operators}
Subspace theorems are concerned with the following phenomenon.
For a well behaved linear transformation  $ T $ on a
Banach space $E$ there exists a  --large-- subspace $F \sbe E $ 
on which $T$ is much better behaved. 
The best known examples of subspace theorems include bounded and non weakly compact operators on the spaces 
$ C(K) ,$ $A,$ $L^\infty,$ $H^\infty,$ and embedding operators on 
$L^p $ spaces.
\begin{enumerate}
\item
If $T$ is a bounded operator on $C(K)$  and  not weakly compact,
then there exists a subspace $F \sbe C(K) $ isomorphic to $c_0$ 
so that $T_{|F}$ is an isomorphism. \cite{pe} The same holds for the disk algebra $A .$
\cite{d}
If T is a  bounded operator on $L^\infty$  and  not weakly compact,
then there exists a subspace $F \sbe  L^\infty$ isomorphic to $\ell^\infty$ 
so that $T_{|F}$ is an isomorphism. \cite{ros} The same assertion holds for operators
on the space $H^\infty. $ \cite{b}
\item
Another class of subspace theorems concern embedding operators  
$ T : L^p \to L^p $ with $ 1 \le p < \infty . $ To any such embedding there exist 
a subspace $ F \sbe   L^p $ so that $F$ and $T(F)$ are  complemented in 
$L^p, $ and $F$ is isomorphic to the 
ambient space  $L^p . $ \cite{es}, \cite{ma}, \cite{jmst}.
Extensions of this theorem hold for rearrangement invariant Banach spaces 
in which the Haar system is an unconditional basis. \cite{jmst}
\end{enumerate}

In this section we  prove a subspace theorem for a 
rearrangement operator acting on the Haar system 
$$ T ( h_I) =  h_{\t(I)} $$ where $\t : \cD \to \cD $ is bijective satisfying
\begin{equation}\numer{10okt21} 
 |\t(I)| = |I| , \quad I \in \cD .
\end{equation} 
\begin{theor}\numer{summ} 
 Assume that  $T$ is bounded on
$L^{p_0}  $ for some  $1 < p_0  \ne2 < \infty .$
Then for any 
$ 1< p < \infty \, (sic!)$ there exists a closed subspace $X_p \sbe
L^p $ isomorphic to $L^p ,$ so that 
$$  T_{|_{X_p}} : (X_p , \| \cdot\|_{L^p}) \to (T(X_p),  \| \cdot\|_{L^p}) $$
is an isomorphism,
$$  \left\| T_{|_{X_p}} \right\| _{L^p}\cdot \left\| T^{-1}_{|_{T(X_p)}}\right\|_{L^p} <\infty, \quad\quad\text{for}
\quad \quad(1< p < \infty) .$$
The subspaces $X_p$  and  $T(X_p) $ are complemented in $L^p.$

\end{theor}
Theorem~\ref{summ} 
is a direct consequence of Theorem~\ref{select} and Theorem~\ref{equivalence}. 
In the course of its proof we use martingale techniques, most notably 
the inequalities of Paley respectively  Stein (in Burkholder's respectively Bourgain's version, \cite{fw}, \cite{m}). 
Thus the method yields extensions to  vector valued rearrangement
operators so that 
$$  \left\|  T_{|_{X_p}}\otimes \Id _E \right\| _{L_E^p}\cdot 
\left\|  T^{-1}_{|_{T(X_p)}}\otimes \Id _E \right\|_{L_E^p} <\infty, \quad\quad\text{for}
\quad \quad(1< p < \infty) ,$$
whenever $E$ satisfies the UMD property. The significance of this remark is
connected with the examples 
in \cite{gmpre} and \cite{pfxmii}
of scalar valued $L^p$ bounded rearrangement operator $T$ 
satisfying
$$  \left\|  T\otimes \Id _E \right\| _{L_E^p} = \infty \quad\text{and}\quad 
\left\|  T^{-1}\otimes \Id _E \right\|_{L_E^p} =\infty.$$

\paragraph{The combinatorial core.}
The following result, 
was   developed 
for the proof that
the spaces $ L^p , 1 \le p < \infty $ are primary. 
We refer to the work of  Enflo and Starbird \cite{es}, 
Johnson, Maurey, Schechtman, and Tzafriri \cite{jmst}, and 
Enflo via Maurey \cite{ma}. 
It is the main combinatorial tool 
by which we
find {\it two tree structures} 
(one in the domain, another in the range of the operator $T$)
that are compatible with the action of  rearrangement  operators.
\begin{prop}\numer{30okt1}
Let $\nu $ be a measure on $[0,1]$ taking values in 
$\{f \in L^1[0,1] : f \ge 0 \}.$
Assume that  
\begin{equation}
\int_0^1 \max_{I \in \cD_N} \nu (I) dt \ge C^{-1} \quad\text{and} \quad  
 \int_0^1 \nu (I) dt \le C |I|, \quad I \in \cD . 
\end{equation}
Then there exist dyadic trees 
$\{ G_I : I \in \cD \} $ and $\{ F_I : I \in \cD \} $
so that 
\begin{equation}\numer{co}
\int_{ G_I} \nu (F_I) dt \ge  \d |I| ,\quad I \in \cD ,  
\end{equation} 
where $ \d = \d(C) . $
\end{prop}
Our definition of a dyadic tree (as given in the preliminaries section) includes the 
requirement that $G_I$ and $F_I$ can be written as  finite unions of dyadic intervals. 
Lemma 9.8 in \cite{jmst} states just that $F_I$ is a finite union of dyadic intervals.
However, the proof of \cite[Lemma 9.8]{jmst} can  
easily be modified to yield that also $G_I$ is a 
finite union of dyadic intervals.
See  \cite{m2}. 
In the applications of Proposition~\ref{30okt1}  
the vector valued measure $\nu $ carries well structured information on 
weak limits of non linear functionals \cite{jmst}, \cite{am}, \cite{m2}.
The non linearities arise by composing the linear operator under 
investigation with the dyadic square function.
 
The hypothesis of Proposition~\ref{30okt1}  are easily verified with the 
following criterion \cite[Chapter 9]{jmst}.
\begin{prop} Let $\nu $ be a measure on $[0,1]$ taking values in 
$\{f \in L^1[0,1] : f \ge 0 \}.$
Assume that 
\begin{equation}\numer{15okt1}
 0 \le \nu (I) \le 1, \quad \quad   \nu ([0,1]) = 1_{[0,1]} , 
\end{equation}
and 
\begin{equation}\numer{15okt2}
\int_0^1 \sqrt{\nu (I)} dt  \le C |I| ,  \quad \quad I \in \cD, 
\end{equation}
then 
\begin{equation}\numer{15okt3}
\int_0^1 \max_{I \in \cD_N} \nu (I) dt \ge C^{2} .
\end{equation}
\end{prop} 
\proof
First observe that by \eqref{15okt1} and the  additivity of the vector measure,
$$ 1_{[0,1]} 
= \left(\sum_{I \in \cD_N}   \nu (I) \right)^{1/2}. $$ 
Next by arithmetic and the  Cauchy-Schwarz inequality, 
\begin{equation}\numer{15okt4}
\begin{aligned}
1 & = \int_0^1 \left(  \sum_{I \in \cD_N}   \nu (I) \right)^{1/2} dt \\
&\le \int_0^1 \left(  \max_{I \in \cD_N} \nu (I)^{1/2} \right)^{1/2}\left(  \sum_{I \in \cD_N}
\nu (I)^{1/2} \right)^{1/2} dt\\
& \le  \left(\int_0^1   \max_{I \in \cD_N} \nu (I) ^{1/2} dt \right)^{1/2} 
 \left(\int_0^1  \sum_{I \in \cD_N}
\nu (I)^{1/2} dt  \right)^{1/2} . 
\end{aligned}
\end{equation}
By \eqref{15okt2}  we get for the second 
term in the   above expression
\begin{equation}\numer{15okt5}
\begin{aligned}
\int_0^1  \sum_{I \in \cD_N}
\nu (I)^{1/2} dt & \le C \sum_{I \in \cD_N} |I|\\
                   &= C.
\end{aligned}
\end{equation}
Combining \eqref{15okt4} with  \eqref{15okt5} and using Hoelder's inequality 
gives 
\begin{equation}\numer{15okt6}
\begin{aligned}
1 & \le C^{1/2} \left(\int_0^1   \max_{I \in \cD_N} \nu (I)^{1/2} dt \right)^{1/2}  \\
&\le  C^{1/2} \left(\int_0^1     \max_{I \in \cD_N} \nu (I) dt\right)^{1/4} 
\end{aligned}
\end{equation}
\endproof

\begin{theor} \numer{select}
To each rearrangement $\t : \cD \to \cD $ satisfying Semenov's condition
\begin{equation}\numer{1feb091}
|\tau (I) | = |I| \quad\text{and}\quad |\t (Q(I))^*| \le C |I| , \quad I \in \cD. 
\end{equation} 
there exist pairwise disjoint collections of dyadic intervals 
$\{ \cH _I : I \in \cD \} $ so that
\begin{enumerate} 
\item The family $ \cH _I$ consists of pairwise disjoint dyadic intervals 
of equal length, that is,   if $J_1 , J_2 \in  \cH _I $ then $|J_1| = |J_2| $ 
and 
$J_1 \cap J_2 = \es .$   
\item The families $\{ \t(\cH _I)^*  : I \in \cD \} $ respectively
$\{ \cH _I^*  : I \in \cD \} $support  dyadic trees. 
\end{enumerate}
\end{theor}
\proof
In the course of selecting the families $\{ \cH _I : I \in \cD \} $
we exploit Proposition~\ref{30okt1}
and take advantage of the fact that we are working with 
rearrangment operators for which \eqref{1feb091} holds.
\paragraph{Selecting two trees. }
We  review the construction of the 
non negative $L^2 $ valued vector measure 
describing the limiting behavior of the operator $T.$
Here we refer to Chapter 9 of \cite{jmst}.  

Fix $ I \in \cD , n \in  \bN .$ Let 
$$ \cD_n (I) = \{ J \in \cD : J \sbe I , |J| = 2^{-n} \} , \quad\text{and}\quad
 d_n (I) = \sum _{J \in  \cD_n (I)} 1_{\t(J)} . $$
Since $  0 \le  d_n (I) \le 1 , $ the sequence 
$\{ d_n (I) : n \in \bN \} $ has a weak cluster point in 
$L^2 . $
By a 
diagonal argument there  exists a subsequence $(n_k)$
so that for  $I \in \cD $ the sequence   $\{d_{n_k} (I)\}$ converges weakly in $L^2 . $
For $ I \in \cD $ define $\nu (I)$ as $L^2-$weak limit,  
$$  
 d_{n_k} (I) \overset{\o}{\longrightarrow}\nu (I) 
$$
Using that $\t$ is bijective and $|\t(I)| = |I| $ it is straightforward to 
observe  that  
\begin{equation}\numer{obs} 
\nu ([0,1]) = 1 , \quad \nu (I) \le 1 , \quad\text{and}\quad \int_0^1 \nu (I) = |I| . 
\end{equation}
Since $\t$ satisfies Semenov's condition, the linear extension of 
$ T ( h_I)=h_{\t (I)} $ 
defines a bounded operator on $H^1 . $
We use the boundedness of $T$  on $H^1  $ 
to prove that 
\begin{equation} \numer{sqrt}
\int_0^1 \sqrt{\nu (I)} \le \|T\|_{H^1} |I| .
\end{equation}
Mazur's theorem  asserts that there exist
$N  \in  \bN,$ $  a _n \in [0,1]  $ so that 
$ \sum_{n \in A_I}\a_n =1$   and 
\begin{equation} \numer{app}
 \|  \sum_{n = 1}^N a_n d_n(I) - \nu (I) \|_{L^2} \le \e^2 |I| .
\end{equation} 
Define next 
$$ k_I = \sum_{n = 1}^ N  a_n^{1/2}\sum _{J \in  \cD_n(I)} h_J . $$
Since $ d_n(I)^2 = d_n(I)  , $ we have the identity 
\begin{equation} \numer{id}
S(Tk_I) = \left( \sum_{n = 1}^N  a_n d_n(I) \right)^{1/2} .
\end{equation} 
By \eqref{app} and \eqref{id}  we get immediately
$$ \int_0^1 \sqrt{\nu (I)} \le (1 + \e ) \int_0^1 S(Tk_I) .$$
Invoking that $\|   k_I\|_{H^1} = |I| $ we obtain that and 
$\| T  k_I\|_{H^1}\le \| T\|_{H^1} |I|$ hence 
\begin{equation}
\begin{aligned}
 \int_0^1 \sqrt{\nu (I)} &\le (1 + \e )\| T\|_{H^1} |I|.
\end{aligned}
\end{equation}
Since $\e >0 $ is arbitrary we obtain \eqref{sqrt} as claimed. 
Combining \eqref{obs} and \eqref{sqrt} yields these estimates
\begin{equation}\numer{input}
\int_0^1\max_{I \in \cD_n} \nu (I) dt \ge  \| T\|_{H^1}^{-2} \quad \quad\text{and}   \quad\quad
 \int_0^1 \nu (I) dt \le |I| . 
\end{equation} 

Since $\nu $ is a finitely additive set function satisfying   $\int_0^1 \nu (I) dt \le |I|,$  
we may extend the mapping 
$$\nu : \cD \to L^2 ( [0,1]) , \quad I \to \nu (I) $$
to an absolutely continuous vector measure  on the 
$\s-$ algebra generated by $\cD $
so that \eqref{input} holds and 
$$
|A_n| \to 0 , \quad \text{ implies } \quad \int_0^1 \nu (A_n ) dt \to 0,  \quad\quad A_n \in \s(\cD) $$
 By Proposition~\ref{30okt1}  
there exist two dyadic trees 
$\{ G_I : I \in \cD \} $ and $\{ F_I : I \in \cD \} $
so that 
\begin{equation}\numer{comb}
\int_{ G_I} \nu (F_I) dt \ge  \| T\|_{H^1}^{-2} |I| ,\quad I \in \cD . 
\end{equation} 
This completes our summary of the selection process in \cite{jmst}. 
\paragraph{Defining $\cH_I $.}
We turn to drawing consequences of 
\eqref{comb}. 
Here we exploit that the operator generating the vector measure 
$ \{ \nu (I) : I \in \cD \}  $
is defined as a rearrangement
$$ T(h_I) = h_{\t(I)} \quad\text{where}\quad  |\t(I)| =  |I| . $$ 

There exists $N(I) \in \bN $ so that  
$G_I$ is a finite union of intervals in $\cD_{N(I)}. $
We test  the weak limit $ \nu (F_I)$
by integrating against the function $1_{G_I} .$
 There exists $M(I) \ge N(I) $ so that 
\begin{equation}\numer{16okt1}
\begin{aligned}
   \int_{ G_I} d_{M(I)} (F_I) dt  &\ge  \frac12  \int_{ G_I} \nu (F_I) dt.
\end{aligned}
 \end{equation}
Define now the collection
\begin{equation}\numer{hi}
\cH_I = \{ J \in \cD_{M(I)}  : J \sbe F_I , \t(J) \sbe G_I \} .
 \end{equation}

\paragraph{Supporting trees. } 
We first verify that $ \{\t(\cH_I)^*  : I \in \cD \}$ supports the tree 
$\{ G_I : I \in \cD \} . $
Recall that by definition of $\cH_I$ we have the inclusion
$$ \t(\cH_I)^* \sbe G_I . $$
It remains for verify the measure estimate: 
Here we use the identity 
 $$\int_{G_I} d_{M(I)} (F_I) = \sum \{ |\t(J)| : J \in \cH_I \}, $$
together with \eqref{comb} and \eqref{16okt1}. This gives
\begin{equation}\numer{me2}
\begin{aligned}
| \t(\cH_I)^*| &= \sum \{ |\t(J)| : J \in \cH_I \}\\
               &  \ge \frac12  \| T\|_{H^1}^{-2} |I| . 
\end{aligned}
\end{equation}
Next we verify  that $\{ \cH_I^* : I \in \cD \}$ supports the tree $\{ F_I : I \in \cD \}.$
Note first the inclusion 
          $$ \cH_I^* \sbe F_I . $$
Next using the identity $ |\t(I)| = |I| $ we get a reduction to 
\eqref{me2} treated above:
\begin{equation}\numer{me1}
\begin{aligned}
| \cH_I^*| &= \sum \{ |J| : J \in \cH_I \}\\
           &  = \sum \{ |\t(J)| : J \in \cH_I \}\\
            & = | \t(\cH_I)^*|\\
              &  \ge \frac12  \| T\|_{H^1}^{-2} |I| . 
\end{aligned}
\end{equation}
\endproof
The collections  $\{ \cH _I : I \in \cD \} $ are chosen so that Stein's Martingale Inequality 
yields the estimates of Theorem~\ref{equivalence} below.  
As a result the scalar valued estimates of Theorem~\ref{equivalence} remain true  when 
the coefficients are taken from a Banach space with the UMD property.

\begin{theor}\numer{equivalence}  Let $\t : \cD \to \cD $ be a rearrangement satisfying 
$|\tau (I) | = |I| , \quad I \in \cD  $ and define
$$ T(h_I) = h_{\tau(I)} . $$
Assume that $\{ \cH _I : I \in \cD \} $ are  pairwise disjoint collections of dyadic intervals so that
\begin{enumerate} 
\item If $J_1 , J_2 \in  \cH _I $ then $|J_1| = |J_2| $ and 
$J_1 \cap J_2 = \es .$   
\item The families $\{ \cH _I^*  : I \in \cD \} $ and 
$\{ \t(\cH _I)^*  : I \in \cD \} $ support dyadic trees of sets. 
\end{enumerate}
Then the blocked system 
\begin{equation}
\widetilde{h_I} = \sum _{ J \in \cH_I } h_J, \quad I \in \cD , 
\end{equation}
satisfies the following estimates: For $1 < p < \infty,$ and any choice of 
$ x_I \in \bR  ,$
\begin{equation} 
\numer{jan10.e1}
c_p \left\| \sum x_I \widetilde{h_I} \right \|_{L^p} \le 
 \left\| \sum x_I T(\widetilde{h_I}) \right \|_{L^p}\le
C_p \left\| \sum x_I \widetilde{h_I} \right \|_{L^p}, 
\end{equation}
\begin{equation} 
\numer{jan10.e2}
c_p \left\| \sum x_I h_I \right \|_{L^p}  \le 
 \left\| \sum x_I \widetilde{h_I} \right \|_{L^p}\le
C_p \left\| \sum x_I h_I \right \|_{L^p}.
\end{equation}
\end{theor}
\proof
Let $\{F_I\}$ be the dyadic tree supported by $\{ \cH _I^*  : I \in \cD \} $
and $\{G_I\}$ be the dyadic tree supported by $\{ \t(\cH _I)^*  : I \in \cD \} $.
Without loss of generality we may assume that the 
tree of sets $\{F_I\}$ is contained in the interval 
$[0,1]$ and that   $\{G_I\}$ is contained in  $[1,2].$
Let $\cE_N $ be the $\s-$algebra generated by 
the sets
       $$ \{F_I\cup G_I ,\quad I \in \cD_N\} .$$
Since $\{F_I\}$ and $\{G_I\}$ are trees and moreover 
$G_J \cap F_I = \es , I, J \in \cD $ it follows that $\{\cE_N , N \in \bN \}$
is an increasing  sequence of  $\s-$algebras.
Let $\bE_N$ denote the conditional expectation operator induced  by 
$\cE_N .$ Note that the following pointwise estimates hold true.
\begin{equation}\numer{crucial}
 1_{G_I} +   |\widetilde{h_I}|
\le C \bE_N(|\widetilde{h_I}|)
\quad\quad I \in \cD_N .  
\end{equation}

Exploiting the estimate \eqref{crucial} together with the unconditionality 
of the Haar system and Bourgain's version of E.M. Stein's martingale
inequality we get the following estimates. Let $x_I \in \bR $ with 
$I \in \cD $ be a sequence with finitely many entries different from zero.

\begin{equation}
\begin{aligned}
\| \sum_{N \in \bN } \sum_{I \in \cD_N} x_I T(\widetilde{h_I})\|_{L^p}
&\le C_p \int_0^1 
\| \sum_{N \in \bN } \sum_{I \in \cD_N} r_I(t) x_I |T(\widetilde{h_I})|\|_{L^p}dt
\\
&\le 
 C_p  \int_0^1\| \sum_{N \in \bN } \sum_{I \in \cD_N} r_I(t) x_I  1_{G_I}\|_{L^p}dt\\
&\le 
 C_p  \int_0^1\| \sum_{N \in \bN } \sum_{I \in \cD_N} r_I(t) x_I  (1_{G_I}+ |\widetilde{h_I}|)|
\|_{L^p}dt\\
&\le 
 C_p  \int_0^1\| \sum_{N \in \bN } \sum_{I \in \cD_N} r_I(t) x_I 
  \bE_N(|\widetilde{h_I}|) 
\|_{L^p}dt
\\
&\le 
 C_p \| \sum_{N \in \bN } \sum_{I \in \cD_N} x_I 
 \widetilde{h_I} 
\|_{L^p}.
\end{aligned}
\end{equation}

Observe also the pointwise estimates. 
\begin{equation}\numer{crucial11}
 1_{F_I} +   |T(\widetilde{h_I})| \le 
 C \bE_N(|T(\widetilde{h_I})|), \quad\quad I \in \cD_N .  
\end{equation}
By \eqref{crucial11} together with the unconditionality 
of the Haar system and Bourgain's version of E.M. Stein's martingale
inequality we get similarly,
\begin{equation}
\begin{aligned}
 \| \sum_{N \in \bN } \sum_{I \in \cD_N} x_I 
 \widetilde{h_I} 
\|_{L^p} &\le 
C_p \int_0^1\| \sum_{N \in \bN } \sum_{I \in \cD_N} r_I(t) x_I  1_{F_I}\|_{L^p}dt\\
&\le 
 C_p  \int_0^1\| \sum_{N \in \bN } \sum_{I \in \cD_N} r_I(t) x_I  (1_{F_I}+  |T(\widetilde{h_I}|)
\|_{L^p}dt\\
&\le 
 C_p \int_0^1 \| \sum_{N \in \bN } \sum_{I \in \cD_N} r_I(t) x_I   \bE_N(|T(\widetilde{h_I})|)
\|_{L^p}dt\\
&\le 
 C_p  \| \sum_{N \in \bN } \sum_{I \in \cD_N}  x_I  T(\widetilde{h_I})
\|_{L^p}.
\end{aligned}
\end{equation}
This proves \eqref{jan10.e1}.

To prove \eqref{jan10.e2}, 
let $\cF_N $ be the $\s-$algebra generated by 
the sets
       $$ \{F_I ,\quad I \in \cD_N \}.$$
Since $\{F_I\}$ is a tree 
$\{\cF_N , N \in \bN \}$
is an increasing  sequence of  $\s-$algebras.
Let $\bF_N$ denote the conditional expectation operator induced  by 
$\cF_N .$ 
The collection $\{ \cH _I^*  : I \in \cD \} $ supports the tree $\{\cF_I, I \in \cD\}$.
Therefore
the following pointwise estimates hold true.
\begin{equation}\numer{crucialF}
 |\widetilde{h_I}| \le   1_{F_I}\le C \bF_N(|\widetilde{h_I}|), \quad\quad I \in \cD_N .  
\end{equation}

\begin{equation}
\begin{aligned}
 \| \sum_{N \in \bN } \sum_{I \in \cD_N} x_I 
 \widetilde{h_I} 
\|_{L^p}dt&\le 
C_p \int_0^1  \| \sum_{N \in \bN } \sum_{I \in \cD_N} r_I(t) x_I  1_{F_I}
\|_{L^p}dt\\
&\le 
 C_p \int_0^1\| \sum_{N \in \bN } \sum_{I \in \cD_N} r_I(t) x_I 
  \bF_N(|\widetilde{h_I}|) 
\|_{L^p}dt
\\
&\le 
 C_p \| \sum_{N \in \bN } \sum_{I \in \cD_N} x_I 
 \widetilde{h_I} 
\|_{L^p}
\end{aligned}
\end{equation}

%
Finally since both since both $\{I\}$ and 
 $\{F_I\}$  are trees of sets we obtain by a 
measure preserving transformation that 
\begin{equation}
\begin{aligned} 
\int_0^1\| \sum_{N \in \bN } \sum_{I \in \cD_N} r_I(t) x_I  1_{I}
\|_{L^p}&dt\le  
C_p\int_0^1\| \sum_{N \in \bN } \sum_{I \in \cD_N} r_I(t) x_I  1_{F_I}
\|_{L^p}dt\\
& \le C_p\int_0^1\| \sum_{N \in \bN } \sum_{I \in \cD_N} r_I(t) x_I  1_{I}
\|_{L^p}dt.
\end{aligned}
\end{equation}
\endproof
The estimates of Theorem~\ref{equivalence} assert that in $L^p , 1 < p < \infty $ the blocked 
system 
is equivalent to the Haar system, hence  the subspace 
$X_p \sbe L^p $ defined to be the $L^p-$closure of 
$\span \{\widetilde{h_I} :   I \in \cD \}  $  is isomorphic to the ambient space $L^p .$ 
And furthermore on the subspace  $X_p$ the rearrangement operator  $T$ acts as an isomorphism. 
\paragraph{Proof of Theorem~\ref{summ} :}
Consider first $1 < p_0 <2 . $ By Semenov's theorem $ \|T\|_{p_0} < \infty $
implies that
$$ |\t (Q(I))^* | \le C |I| , \quad I \in \cD . $$
Apply Theorem~\ref{select} to the rearrangement $ \t : \cD \to \cD . $
Let $ \{ \cH_I : I \in \cD \} $
denote the collections of dyadic intervals satisfying the conclusion of 
Theorem~\ref{select}. Put 
$$ X = \span\{  \widetilde{h_I} : I \in \cD \} \quad\text{where}\quad 
\widetilde{h_I} = \sum_{ J \in \cH_I } h_J . $$
Next fix $ 1 < p < \infty . $ Let $X_p $ denote the closure of $X$ in $L^p , $ and 
similarly let $(T(X))_p$ be the closure of $T(X)$ in $L^p . $ 
Theorem~\ref{equivalence} asserts that $X_p $ is isomorphic to $L^p ,$ and that the map 
$$ T (\widetilde{h_I}) =  \sum_{ J \in \cH_I } h_{\tau (J)}  $$
extends uniquely to an isomorphism with domain $X_p$ and range $(T(X))_p .$ Denoting the 
extension still by  $T$ we have $(T(X))_p  = T(X_p),$ and 
$$  \left\| T_{|_{X_p}} \right\| _{L^p}\cdot \left\| T^{-1}_{|_{T(X_p)}}\right\|_{L^p} <\infty . $$
 Next we turn to the case $ 2 < p_0 < \infty . $ Note that for rearrangement operators 
the transposed operator coincides with the inverse $T^{-1} $ defined by 
$ \t^{-1} : \cD \to \cD .$ Thus $ S = T^*  = T^{-1}  $ is a bounded operator on 
$L^{p_0^*} $ where $ 1/p_0^* + 1/p_0 = 1 . $ Now  $ p_0^* < 2 . $ 
By the first part of the theorem  applied to the rearrangement operator $S,$ for 
$1 < p < \infty $  
there exists $X_p$ isomorphic to $L^p, $ so that 
$$  \left\| S_{|_{X_p}} \right\| _{L^p}\cdot \left\| S^{-1}_{|_{S(X_p)}}\right\|_{L^p} <\infty . $$
Since $S = T^{-1} ,  S^{-1} = T $ with  $ Y_p = S(X_p)$
we obtain 
$$  \left\| T_{|_{Y_p}} \right\| _{L^p}\cdot \left\| T^{-1}_{|_{T(Y_p)}}\right\|_{L^p} <\infty . $$
\endproof

\section{\numer{sem.s1} Shift operators}

We consider rearrangements of $\cd$ of particular type, namely shift operators. Let 
$$
I_{j,k} = [{k-1 \over 2^j}, {k \over 2^j}) \quad {\rm for} \quad k=1, \ldots, 2^j,
$$
and let
$M=\{m_j, j \geq 1\}$ be a sequence of integers
satisfying $|m_j| \leq 2^{j}$. Consider $\tau_M : \cd \to \cd$ given by 
$$\tau_M(I_{j,k}) = I_{j,k+m_j} = I +m_j|I| \quad {\rm for} \quad
I \in \cd_j,
$$
where $k+m_j$ is understood ${\rm mod} \; 2^j$. This rearrangement is called {\em shift  generated by $M$}.


 

In this section, we give a version of Semenov's theorem suitably adapted to the special nature of shift operators, see 
Proposition \ref{sem.p1}. We isolate a class of shift operators for which Conjecture \ref{con1} and Conjecture \ref{con2} hold true,
Theorems \ref{sem.t2} and  \ref{sem.t1}. We prove directly, without using the machinery of Section \ref{sub} 
a subspace theorem for shift operators, Theorem \ref{sem.t3}.

\subsection{Semenov condition for shifts.}
We next give a version of Semenov's criterion that holds specifically for shift operators.
Let $j \in \bN$ and  put 
$$ x_j = {m_j \over 2^j},
$$
so that 
$$\tau_M(I) = I + x_j \quad {\rm for} \quad |I|= 2^{-j}.
$$
Then define
$$
N_j(M) = | \{k: I_{j,k} \cap \{x_l, l\geq j\} \neq \emptyset \} |.
$$
The next Proposition \ref{sem.p1} relates Semenov's condition to the boundedness of the sequence $N_j(M)$.
\begin{prop}
\numer{sem.p1}
Let $\tau_M$ be the shift generated by the sequence $M=\{m_j, j \geq 1\}$. 
Then we have for $I \in \cd_j$
\begin{equation}
\numer{w.e1}
{N_j (M)\over 2} |I| \leq |\tau_M(Q(I))^*| \leq 2N_j(M) |I|.
\end{equation}
Consequently,  $\tau_M$ satisfies Semenov's condition
if and only if there is a constant $K>0$ such that $N_j(M) \leq K$ for all $j \geq 1$.
\end{prop}
\proof
Let  $N_j = N_j(M)$. For fixed $j$, let  $$1 \leq k_1 <k_2 < \ldots < k_{N_j} \leq 2^j$$
be the enumeration of the set of indices
$\{k: I_{j,k} \cap \{x_l, l\geq j\} \neq \emptyset \}$.
Take $I \in \cd_j$ and
$l \geq j$. Let $k_i$ be such that $x_l \in I_{j,k_i}$.
Then we have 
$$
\tau_M(Q(I) \cap \cd_l) ^* = I + x_l   \subset (I + {k_i-1 \over 2^j}) \cup (I + {k_i \over 2^j}),
$$
and consequently
$$
| \tau_M(Q(I))^*| \leq 2N_j |I|.
$$
To prove the other estimate, let $J,L \in \{I_{j,k_i}, 1 \leq i \leq N_j\}$ be such that
${\rm dist}(J,L) >0$. Then for $l_1, l_1$ such that $x_{l_1} \in J$ and $x_{l_2} \in L$ we have
$$
\tau_M(Q(I) \cap \cd_{l_1}) ^* \cap \tau_M(Q(I) \cap \cd_{l_2}) ^* = \emptyset.
$$


Clearly, for each $l \geq j$ we have $|\tau (Q(I) \cap \cd_l)^*| = |I|$. 
Combining these observations we get

$$
| \tau_M(Q(I))^* | 
\geq {N_j \over 2} |I|.
$$
\endproof



\subsection{Shifts and nested collections.}
Proposition \ref{sem.p1} implies that if $\tau_M$ satisfies Semenov's condition, then the collection of accumulation points of
sequence $X = \{x_j, j \geq 1\}$ is finite. 
Therefore, without loss of generality we assume $\lim_{j \to \infty} x_j = 0$.

We turn our attention to sequences $M$ with $\limsup_{j \to \infty} N_j(M) \leq 2$.
We are able to describe the structure of such sequences:

\begin{prop}
\numer{sem.p2}
Let  $M=\{m_j, j \geq 1\}$ and $X= \{x_j= {m_j \over 2^j}, j \geq 1\}$ 
be such that $\lim_{j \to \infty} x_j =0$ and $\limsup_{j\to \infty} N_j =  2$.
Then exist  sequences $$\{a_k, k \geq 1\}
\quad {\it and} \quad \{j_k, k \geq 0\} $$ 
with the following properties: $\{j_k, k \geq 0\}$ is increasing 
and
if $$j_{k-1} \leq j < j_{k}$$ then
$$
x_j = 0 \quad {\it or} \quad
 | a_k - x_j| < {1 \over 2^j},
$$
and if 
$$j \ge j_k $$
then
$$
x_j \le 2^{-{j_k}+1} .
$$
\end{prop}

\proof 
First, note that if $N_j =1$ then $0 \leq x_l < {1 \over 2^j}$ for all $l \geq j$.
Therefore, in case $\limsup_{j \to \infty} N_j =1$ the above condition is clearly satisfied.

Consider the case $\limsup_{j \to \infty} N_j =2$. 
We construct  sequences $\{a_k, k \geq 1\}$ and $\{j_k, k \geq 0\}$ inductively. We put an additional requirement that $N_{j_k}=2$.
Let
\begin{eqnarray*}
j_0 & = & \min\{j: N_i \leq 2 \quad \hbox{for all} \quad i \geq j \quad \hbox{and} \quad N_j =2 \},
\\
n_1 & = & \max\{j \geq j_0: x_j > {1 \over 2^{j_0}}\},
\\
j_1 & = & \min\{j > n_1: N_j = 2\}.
\end{eqnarray*}
For  $n_1 < j < j_1$ we have $N_j = 1$, and consequently $\{ x_l, l \geq j\} \subset [0, {1 \over 2^j})$.
Taking $j=j_1-1$ we get $x_l \leq {1 \over 2^{j_1-1}}$ for $l \geq j_1$.

Then, take $l$ such that $j_0 \leq l \leq n_1$. Since $N_l =2$, there are only two intervals
from $\cd_l$ which have nonempty intersection with $\{ x_i, i \geq l\}$: one of them is $I_{l,1}$, and the other 
contains $x_{n_1}$. This implies that either $x_l \in I_{l,1}$ (in such case, $x_l=0$), or $x_l$ and $x_{n_1}$
must be in the same interval from $\cd_l$, which implies $|x_{n_1} - x_l| < {1 \over 2^l}$.
Therefore, it is enough to put $a_1 = x_{n_1}$.
Note that $N_{j_1} =2$.

Having defined $j_k$ and $a_k$, we define $j_{k+1}$ and $a_{k+1}$. Since $N_{j_k} =2$, there is $j\geq j_k$ such that
$x_j > {1 \over 2^{j_k}}$. 
We put
\begin{eqnarray*}
n_{k+1} & = & \max\{j \geq j_k: x_j > {1 \over 2^{j_k}}\},
\\
j_{k+1} & = & \min\{j > n_{k+1}: N_j = 2\},
\\
a_{k+1} & = & x_{n_{k+1}}.
\end{eqnarray*}
Arguments analogous to the above one show that Proposition \ref{sem.p2} 
holds for this choice of $j_{k}$, $j_{k+1}$ and $a_{k+1}$.
\endproof

\begin{theor}
\numer{sem.t2}
Let $M$ be a sequence such that $\limsup_{j\to \infty} N_j(M) \leq 2$.
Then the operator $T_M h_I = h_{\tau_M(I)}$ extends to an isomorphism of 
$L^p_E$, for each Banach space $E$ with UMD property and $1<p<\infty$.
\end{theor}

\proof Without loss of generality, we assume that on each $\cd_j$ the shift is nontrivial,
so that $m_j \neq 0$ and consequently $x_j \neq 0$, $j \geq 1$.
By Proposition \ref{sem.p2}, the hypothesis  
$$
\limsup_{j\to \infty} N_j(M) \leq 2
$$
implies that there exist
$$\{a_k, k \geq 1\}
\quad {\rm and} \quad \{j_k, k \geq 0\}, $$ 
where $j_k, k \geq 0$ is increasing and satisfies the following properties.
\begin{equation}
\numer{cc.e2}
\text{ If } j_{k-1} \leq j < j_{k} \text{ then }
 | a_k - x_j| < {1 \over 2^j},
\end{equation}
and
\begin{equation}
\numer{cc.e3} 
j \ge j_k  \text{ implies }
x_j \le 2^{-{j_k}+1} .
\end{equation}

Applying   an additional shift  by 
no more than $1$ unit, we can assume that for each $j_{k-1} \leq j < j_k$ and $I \in \cd_j$
$$
|(x_j + I) \cap (a_k +I) |\geq {|I| \over 2}, 
$$
while for $j \geq j_k$ we have $x_j \leq {3 \over 2^{j_k}}$. 
Recall that by Theorem \ref{tf.t1}, the additional shift applied above induces an isomorphism of $L^p_E$,
with uniform bounds.

The idea is to 
take a suitable splitting of level $\cd_{j_{k-1}}$, and for $\cd_j$ with $j_{k-1} \leq j <j_k$
take a splitting induced by the splitting of $\cd_{j_{k-1}}$: two intervals from $\cd_j$ will be put into the same collection 
in the splitting of $\cd_j$
iff their dyadic ancestors in $\cd_{j_{k-1}}$ are put into the same collection in the splitting of $\cd_{j_{k-1}}$. Then
we would like to
take $I \cup (a_k +I)$ as our building block for the elements in the nested collection for $I \in \cd_j$ with $j_{k-1} \leq j <j_k$.
In order to make this work, we need to produce space between the sets $I \cup (a_k +I)$.
Therefore
we perform now one additional preparatory operation.
Define $\phi:\cd \to \cd$ by the following procedure: given $I \in \cd$, let
$r(I)$ be the right endpoint of $I$ and $s(I)$ be its midpoint. 
Then define $\phi(I)$ uniquely by the following relation:
$$
|\phi(I)| = {1 \over 4} |I| \quad \text{and} 
\quad r\big(\phi(I)\big) = s(I).
$$
That is, $\phi:\cd \to \cd$ be the mapping assigning to  $I \in \cd_j$ an interval $\phi(I) \in \cd_{j+2}$
with right endpoint coinciding with the midpoint of $I$.
The reason we use the mapping $\phi$ is that the intervals in the range of $\phi$ satisfy the following property:
\begin{equation}
\numer{cc.e1}
\text{ If } J_1, J_2 \in \phi(\cd) 
\text{ with } |J_1| \leq |J_2|, \text{ then } J_1 \subset J_2
\text{ or } {\rm dist}(J_1,J_2) \geq |J_1|.
\end{equation}

Let $\psi:\cd \to \cd$ be the mapping assigning to $I \in \cd_j$ the unique interval $\psi(I) \in \cd_{j-2}$
such that $I \subset \psi(I)$. 
The operations $T_\phi$ and $T_\psi$ are bounded on $L^p_E$. 
Let 
\begin{equation}
\numer{jj}
j_k'=j_k+2.
\end{equation}
Let $\sigma$ be a shift such that for $I \in \cd_{j}$ with $j_{k-1}' \leq j < j_k'$ 
$$
|\sigma (I) \cap (a_k +I)| \geq {|I| \over 2}.
$$
Put $\tau' = \psi \circ \sigma \circ \phi$. Note that $\tau$ can be obtained from $\tau'$ by applying  
an additional shift by at most 1.  
Therefore, to estimate the norm of 
$T_\tau$, it is enough to consider $\sigma$ restricted to $\phi(\cd)$.
In the next paragraph, we will split $\phi(\cd)$ as 
$$
\phi(\cd) = \bigcup_{r=1}^6 \bigcup_{l=1}^{512} \cf(r,l)
$$
so that the restricted rearrangements $\sigma:\cf(r,l) \to \cd$ are supporting nested collections.
Combining this with Theorem \ref{sem.p20}, proves the statement of 
Theorem \ref{sem.t2}.

\paragraph{Horizontal splitting at stage $k.$}
In this paragraph we fix $k.$ We will  obtain a splitting 
of the collection 
$$ \phi\left(\bigcup_{j=j_{k-1}}^{  j_k-1 } \cD_j\right) .$$
We obtain it by  first decomposing the collection of 
top level intervals $\phi(\cD_{j_{k-1}})$ and then simply pushing it down to later levels.
(Hence the name horizontal splitting.) 
On $\cd_{j_{k-1}}$, $\tau$ is a shift by $m_{j_{k-1}}$. We  split $\cd_{j_{k-1}}$ as 
$$\cd_{j_{k-1}} = \bigcup_{l=1}^L \cd_{j_{k-1},l}, \quad L \leq 8^3$$
so that if 
$$I \in \cd_{j_{k-1},l}$$ 
then 
\begin{equation}
\numer{cc.e5}
I \pm |I| \not \in \cd_{j_{k-1},l}, 
\quad
\tau(I), \tau(I) \pm |I|, \not \in \cd_{j_{k-1},l}, \quad \text{and} \quad \tau(I) \pm 2|I| \not \in \cd_{j_{k-1},l}.
\end{equation}
This is done in a straightforward manner by consecutive separation along the orbits of the three  shifts  $m_{j_{k-1}}$ and 
$m_{j_{k-1}} \pm 2$.
The number of the collections obtained this way admits an universal bound $L \le 8^3.$

We next employ the collections $\cd_{j_{k-1},l}$ to split each of collections $\cd_j \cap \phi(\cd)$ with
 $j_{k-1}' \leq j < j'_k$. That is, we define $\cb_{j,l}$ as the collection of intervals in $\cd_j \cap \phi(\cd)$
whose dyadic predecessor in $\cd_{j_{k-1}}$ is actually contained in $\cd_{j_{k-1},l}$. Thus
\begin{equation}\numer{cc.e6}
\cb_{j,l} = \{ J \in \cd_j \cap \phi(\cd): J \subset \cd_{j_{k-1},l}^*\}.
\end{equation}

Now we fix $l \le L $ and $k$ as above together with $j_k' = j_k +2$ and $j_{k-1}' = j_{k-1} +2,$
to analyze the joint properties of the intervals in the collection
$$\cG = \cup_{j=j_{k-1}'}^{j_k'-1} \cb_{j,l}.$$
If $I,J \in \cG$ with $|I| \leq |J|$, then by \eqref{cc.e1} we have:
$$\text{Either }  I \subset J \text{  or } {\rm dist}(I,J )\geq |I|.$$
Clearly by shifting $I,J$ with $a_k$ this implies that 
 either $a_k+I \subset a_k+J$, or
${\rm dist}(a_{k}+I, a_k+J) \geq |I|.$
Next we exploit the condition \eqref{cc.e5}. 
Together with \eqref{cc.e6} condition \eqref{cc.e5} gives
$$
{\rm dist}(I, a_k +J)\geq {1 \over 2^{j_{k-1}}}\quad\text{and}\quad {\rm dist}(a_k+I, J)  \geq {1 \over 2^{j_{k-1}}}.
$$
Recall that ${1 \over {2^{j_{k-1}'}}}$ is in fact
the length of the largest interval in $\cG$.
Therefore, the distance between $I$ and $J+a_k$ is always larger than $4$ times the length of the largest interval in $\cG$.
Hence if $I \cup (I+a_k)$ is not contained in $J \cup (J+a_k)$, then
$I \cup (I+a_k)$ and $J \cup (J+a_k)$ are separated by more than $\min \{|I|,|J|\}.$
Summing up we arrived at  the following alternative for $I,J \in \cG $ with $|I|\le |J|$ then: 
\begin{equation}
\numer{tf.e14}
\text{Either }\quad I \cup( I + a_k) \subset J \cup (  J +a_k ) \quad\text{ or }\quad {\rm dist}\big(I \cup (a_k+I), J \cup (a_k+J)
\big)\geq |I|
\end{equation}
Next we compare the above separation condition with the diameter of the smallest interval containing $I$ and
$I+a_k$. By \eqref{cc.e2} we have
\begin{equation}
\numer{tf.e10}
{\rm diam} (I \cup a_k+I) \leq {m_{j_{k-1}+}+1 \over 2^{j_{k-1}}} = x_{j_{k-1}} + {1 \over 2^{j_{k-1}}}.
\end{equation}
Later in the proof, we will exploit that the separation \eqref{tf.e14} at stage $k$ is much wider
than the diameter in \eqref{tf.e10} at stages $k+5$ and following. It is only this  implication which makes our construction work.
It is here where we rely on the strong dichotomy expressed by our hypothesis \eqref{cc.e2}-\eqref{cc.e3}. 
Above $k\in \bN $ and $l \le L $ were fixed. We write now 
$$ \cG (k,l) = \cG .$$
Thus we obtained the decomposition 
\begin{equation}
\numer{jan09.e1}
\phi\left(\bigcup_{j=j_{k-1}}^{  j_k-1 } \cD_j\right) =\bigcup_{l=1}^L \cG (k,l), \quad \quad L \leq 8^3.
\end{equation}

\paragraph{Consequences of \eqref{tf.e10}.} 
Here we specify the form of the diameter estimates \eqref{tf.e10} at stages later than $k.$
Fix as before $l \le L. $ Then for
$$h > k \quad\text{ form }
\quad \cG (h,l )= \cup_{j=j_{h-1}'}^{j_h'-1} \cb_{j,l}.$$
Note that if  
$$j_{h-1}' \leq j < j_h' \quad \text{then}\quad  x_j \leq {3 \over 2^{j_{h-1}}} .$$
Therefore by \rf{tf.e10} if $I \in \cG (h,l )$ we get 
$$
{\rm diam} (I \cup a_h+I)  \le \left[{4 \over 2^{h-k-2}}\right] {2 \over 2^{j'_{k}}}.
$$
Next observe that if $ h $ is  larger than $k+6$ then the left hand factor on the right hand side of the above estimate 
is bounded by $1/4 .$
Hence for $h \ge k+6$ and $I \in \cG (h,l )$ we get 
\begin{equation}
\numer{tf.e11}
{\rm diam} (I \cup a_h+I) \leq {1 \over 4} {2 \over 2^{j'_{k}}}.
\end{equation}
\paragraph{Vertical splitting.} It follows from \eqref{jan09.e1} that
$$
\phi(\cd) = \bigcup_{l=1}^L \bigcup_{k=0}^\infty \cG(k,l).
$$
Now, we split the sequence $\{j_k, k \geq 1\}$
into $6$ subsequences $\{j_{r+6s}, s \geq 0\}$, $r=1,\dots, 6$. 
We put
$$
\cF_s(r,l) = \cG(r+6s,l),
$$
thus
$$
\phi(\cd) = \bigcup_{l=1}^L \bigcup_{r=1}^6 \bigcup_{s=0}^\infty \cF_s(r,l).
$$

\paragraph{Construction of a tree.} Fix $1 \leq r \leq 6$ and  $1 \leq l \leq L$.
Let 
$$\cf_s = \cF_s(r,l)\quad\text{ and }\quad k_s = r+ 6s .$$   
For $I \in \cf_s$ and $n\geq s$, we construct two  sets $B_n(I)$ and $ C_n(I)$,  with the following
properties:
\begin{equation}
\numer{11jan09.e1}
I  \subset  B_n(I)  \quad \text{ and }  \quad a_{k_s}+I  \subset  C_n(I) ,
\end{equation}
and inversely 
\begin{equation}
\numer{11jan09.e2}
B_n(I) \subset  \big\{t: {\rm dist}(t, I) \leq  F_{n-s} |I|\big\} \quad
\text{ and }  \quad
C_n(I) \subset  \big\{t: {\rm dist}(t, a_{k_s}+I) \leq F_{n-s} |I|\big\},
\end{equation}
where $F_{n-s} = 2 \sum_{i=1}^{n-s} 4^{-i} $. Thus $B_n(I) $ is contained in an interval with the same midpoint
as $I$ and diameter bounded by $(1+ 2 F_{n-s})\cdot |I| ,$ and the same for
$C_n(I) $ and $a_{k_s}+I.$  This will give us a nested collection defined by 
                $$ A_n(I) = B_n(I) \cup C_n(I) ,\quad I \in \cF_0\cup\dots\cup \cF_n$$ 
so that for 
 $I,J \in \cF_0\cup\dots\cup \cF_n$  and  $|I| \leq |J|$ 
we have:
\begin{equation}
\numer{11jan09.e3}
\text{Either }\quad A_n(I) \subset A_n(J) \quad\text{ or }\quad {\rm dist}(A_n(I), A_n(J)) \geq {2 \over 2^{j_{k_n}'}}.
\end{equation}
We recall that the factor ${2 \over 2^{j_{k_n}'}}$ is the length of the shortest  intervals in $\cF_n.$ 

The construction is inductive. 
\paragraph{Step $n=0$.} For $I \in \cf_0$, put
$$
B_0(I) =  I, \quad C_0(I) = a_{k_0} +I, \quad A_0(I) = B_0(I) \cup C_0(I).
$$
\paragraph{\bf Step $n+1$.}
In step $n$, we defined  $ B_n(I), C_n(I)$, for $I \in \bigcup_{s=0}^n \cf_s .$
For $J \in \cf_{n+1}$ we put
$$
B_{n+1}(J) = J, \quad C_{n+1}(J) = a_{k_{n+1}} + J, \quad A_{n+1} (J) = B_{n+1}(J) \cup C_{n+1}(J).
$$
 Now we define  $B_{n+1}(I) $ and $ C_{n+1}(I)$ for $ I \in  \bigcup_{s=0}^n \cf_s $ by updating
  $ B_n(I), C_n(I).$ To this end we define the index sets 

  $$
  \begin{aligned}
K_{n+1}(I) & =  \{ J \in \cf_{n+1}: {\rm dist}(B_n(I), A_{n+1}(J)) \leq { 2 \over 2^{j_{k_{n+1}}'}}\},
\\
L_{n+1}(I) & = \{ J \in \cf_{n+1}: {\rm dist}(C_n(I), A_{n+1}(J) ) \leq { 2 \over 2^{j_{k_{n+1}}'}}\},
\end{aligned}
$$
and do the updating
$$
  \begin{aligned}
B_{n+1}(I) & = B_n(I) \cup \bigcup_{J \in K_{n+1}(I)}  A_{n+1}(J),
\\
C_{n+1}(I) & = C_n(I) \cup  \bigcup_{J \in L_{n+1}(I)}  A_{n+1}(J),
\end{aligned}
$$
We complete the definition of the tree by putting
$$
A(I) = \bigcup_{n=s}^\infty A_n(I) \quad {\rm for} \quad I \in \cf_s.
$$
It follows by \eqref{11jan09.e3} that if $I,J \in \cf$ and
$|I| \leq |J|$, then either $A(I) \subset A(J)$ or $A(I) \cap A(J) = \emptyset$.
Moreover, as a consequence of \eqref{11jan09.e1}-\eqref{11jan09.e2} we have
$$
I \subset A(I), \quad |\sigma(I) \cap A(I)| \geq {1 \over 2} |I|
\quad {\rm and} \quad  2|I| \leq |A(I)| \leq {20 \over 3} |I|.
$$

\paragraph{Verification of \eqref{11jan09.e1}-\eqref{11jan09.e3}.}
It remains to check \eqref{11jan09.e1}-\eqref{11jan09.e3}. The proof is inductive.

For $n=0$ and $I,J \in \cf_0$, \eqref{11jan09.e3} follows by \eqref{tf.e14}, while
\eqref{11jan09.e1}-\eqref{11jan09.e2} are immediate consequences of the defintion of $B_0(I), C_0(I)$.

Assume that \eqref{11jan09.e1}-\eqref{11jan09.e3} hold at stage $n$. To verify them at stage $n+1$,
recall that for $J \in \cf_{n+1}$ we have
${\rm diam} A_{n+1}(J) \leq {1 \over 4} {2 \over 2^{j'_{k_n}}}$, cf. \eqref{tf.e11}. 
For $I \in \cf_s$ we have $|I| \geq {2 \over 2^{j'_{k_s}}}$. Since $k_n = k_s + 6(n-s)$
it follows that $j'_{k_n} \geq 6(n-s) + j'_{k_s}$.
Therefore, using the induction hypothesis
on $B_n(I)$ we get for $t \in A_{n+1}(J)$ with $J \in K_{n+1}(I)$ 
$$
{\rm dist}(t,I) \leq {\rm diam} A_{n+1}(J)  + { 2 \over 2^{j_{k_{n+1}}'}} + 2
\sum_{i=1}^{n-s} 4^{-i} |I| \leq 2 \sum_{i=1}^{n+1-s} 4^{-i} |I|.
$$
Condition (i) for $C_{n+1}(I)$ is checked in the same way.

To check condition (ii), we consider several cases:
\begin{itemize}
\item[(a)] If  $I,J \in \cf_{n+1}$ and $I \neq J$, then  
${\rm dist}(A_{n+1}(I), A_{n+1}(J)) \geq {2 \over 2^{j'_{k_{n+1}}}}$ by \rf{tf.e14}.

\item[(b)] If $I \in \cf_s$ with $s \leq n$ and $J \in \cf_{n+1}$ is such that $J \in K_{n+1}(I) \cup L_{n+1}(I)$,
then $A_{n+1}(J) \subset A_{n+1}(I)$, by definition of $A_{n+1}(I)$.

\item[(c)] If $I \in \cf_s$ with $s \leq n$ and $J \in \cf_{n+1}$ is such that $J \not\in K_{n+1}(I) \cup L_{n+1}(I)$,
then ${\rm dist} (A_n(I), A_{n+1}(J)) \geq {2 \over 2^{j'_{k_{n+1}}}}$
by the definition of $K_{n+1}(I)$ and $L_{n+1}(I)$.
If $J' \in K_{n+1}(I) \cup L_{n+1}(I)$, then
${\rm dist} (A_{n+1}(J), A_{n+1}(J')) \geq {2 \over 2^{j'_{k_{n+1}}}}$ by \rf{tf.e14}.
Therefore $${\rm dist}(A_{n+1}(I), A_{n+1}(J)) \geq {2 \over 2^{j'_{k_{n+1}}}}.$$

\item[(d)] If $I,J \in \bigcup_{s=0}^n \cf_s$ are such that $A_n(I) \subset A_n(J)$ then $K_{n+1}(I) \cup L_{n+1}(I)
\subset K_{n+1}(J) \cup L_{n+1}(J)$, and consequently $A_{n+1}(I) \subset A_{n+1}(J)$.

\item[(e)] Finally, let $I,J \in \bigcup_{s=0}^n \cf_s$ be such that ${\rm dist}(A_n(I), A_n(J)) \geq {2 \over 2^{j'_{k_n}}}$.
Recall that for $J' \in K_{n+1}(I) \cup L_{n+1}(I)$ or $J' \in K_{n+1}(J) \cup L_{n+1}(J)$ we have
${\rm diam} A_{n+1}(J') \leq  {1 \over 4} {2 \over 2^{j'_{k_n}}}$. Since $k_{n+1} = k_n + 6$
we have $j'_{k_{n+1}} \geq  j'_{k_n} + 6$, and we get 
\begin{eqnarray*}
{\rm dist}(A_{n+1}(I), A_{n+1}(J))& \geq & {\rm dist}(A_n(I), A_n(J)) - 2 ({2 \over 2^{j'_{k_{n+1}}}}
+  {1 \over 4} {2 \over 2^{j'_{k_n}}} )
\\ & \geq & 
{1 \over 2^{j'_{k_n}}} - {4 \over 2^{j'_{k_{n+1}}}}
\geq {2 \over 2^{j'_{k_{n+1}}}}.
\end{eqnarray*}
\end{itemize}

This completes the construction of a tree for $\sigma$ restricted to $\cf$.\endproof


\paragraph{Subspace theorem for shifts.}
As immediate application of the above theorem we prove now 
that the subspace theorem holds for shift operators in a very peculiar way.
Without assuming that the shift operator itself is bounded on $L^p$ we are able to 
find a subspace of $L^p$ on which  the operator acts as an isomorphism.
Moreover our argument here does not use  any of the construction developed 
in Section \ref{sub}.
\begin{theor}
\numer{sem.t3}
Let $M = \{m_j, j \geq 1\}$ be a sequence of integers satisfying $|m_j| \leq 2^j$,
and let $\tau_M$ be the associated shift on $\cd$.
Then there is a sequence $\{j_k, k \geq 1\}$ such that 
$T_M h_I = h_{\tau_M(I)}$ extends to an isomorphism of 
${\rm span}_{L^p_E} \{x_I h_I, I \in \bigcup_{k=1}^\infty \cd_{j_k}, X_I \in E\}$
for each Banach space $E$ with UMD property and $1<p<\infty$.
\end{theor}

\proof Observe that the sequence $x_n = {m_n \over 2^n}$ has an accumulation point.
Without loss of generality we can assume that it has a subsequence $n_l$
such that 
$\lim_{l \to \infty} x_{n_l} = 0$ and $x_{n_l} \neq 0$.
The sequence $(j_k, k \geq 1)$ is defined inductively:
$j_1 = n_1$. When $j_1, \ldots, j_k$ are already defined, we put
$$
j_{k+1} = \min\{n_l > j_k: 0 \leq x_{n_s} < {1 \over 2^{j_k}} \hbox{ for all } s \geq l\}.
$$
Consider sequence $M' = \{m_j', j \geq 1\}$ defined by $m'_{j} = m_j$ if $j=j_k$ for some $k$
and $m_j'=0$ otherwise.
By the definition of the sequence $j_k$ we have
$N_{j_k}(M') = 2$ and $N_j(M') =1  $ for $j \neq j_k$.
Application of Theorem \ref{sem.t2} completes the proof.\endproof


\paragraph{Decomposable sequences.}
Let us fix a shift operator given by  $M=\{m_j, j \geq 1\}.$ 
By the following condition  we attempt to capture the essence of the conclusion 
in  Proposition \ref{sem.p2} and at the same time we would like to allow for a higher degree of flexibility. 
The condition is chosen so that the proof given for Theorem~\ref{sem.t2} actually shows 
that shifts satisfying the condition below induce isomorphisms on $L^p_E$, for each  UMD-space $E$ 
(See Theorem~\ref{sem.t1}.)
\begin{defi}
\numer{sem.d1}
Let  $M=\{m_j, j \geq 1\}$ and $X= \{x_j= {m_j \over 2^j}, j \geq 1\}$ 
be such that $\lim_{j \to \infty} x_j =0$. We say that the sequence $M$ is  decomposable if there are $w_1, w_2 \geq 1$,
 a sequence $\{a_k, k \geq 1\}$ and an increasing sequence of natural numbers
$ \{j_k, k \geq 0\}$ such that
$$
x_j \le {w_1 \over 2^j} \quad {\it or} \quad
 | a_k - x_j| < {w_2 \over 2^j} \quad {\it for } \quad j_{k-1} \leq j < j_{k} ,
$$
and 
$$
x_j \le {w_1 \over  2^{{j_k}-1}} \quad {\it for} \quad j \geq j_k.
$$
\end{defi}
As stated above, the argument given in the course of proving Theorem \ref{sem.t2}
can be adapted in a straightforward way to provide the proof of the following result.
\begin{theor}
\numer{sem.t1}
Let $M$ be a decomposable sequence. 
Then the operator $T_M h_I = h_{\tau_M(I)}$ extends to an isomorphism of 
$L^p_E$, for each Banach space $E$ with UMD property and $1<p<\infty$.
\end{theor}
Proposition \ref{sem.p2} states that a sequence $M$ with 
$\limsup_{j \to \infty} N_j(M) \leq 2$ is decomposable in the sense that it satisfies the condition of Definition~\ref{sem.d1}.
To close this section, we fomulate the following conjecture:

\begin{conjecture}
Each  shift satisfying 
the Semenov-type condition 
$$\sup N_j(M) < \infty $$
can be written as a finite composition of shifts 
satisfying the  hypothesis of  Definition~\ref{sem.d1}.
\end{conjecture}

By Theorem \ref{sem.t1}, the positive answer to this conjecture implies the positive answer to Conjecture \ref{con1} for shift operators.

\section{\numer{c.s1} Combinatorics of
coloured intervals}
The intricacies of the extension problem as formulated in Conjecturies \ref{con1} and \ref{con2}
give rise to the following two-person game of general combinatorial interest.
The game is played by two players with collections of coloured dyadic intervals in 
$$\cD_j= \{ I \in \cD : |I| = 2^{-j} \} \quad  
\text{for a fixed  $j \in \bN. $}$$
It starts by fixing $\eta > 0,$  $d\in \bN ,$ and a subcollection 
$$
\cC(0) \subset \cD_j
$$
with with an $(\eta,d)$-homogeneous colouring
$$\cc_1(0), \ldots , \cc_d(0).$$
(see Definition \ref{c.d1} below). The rules of the game are as folows: Throughout the game, $j \in \bN$ is fixed.
\begin{enumerate}
\item In the first stage, Player A chooses a collection $\cC(1) \varsupsetneq \cC(0)$ and $\cC(1) \subset \cD_j$.
Player B determines an $(\eta,d)$-homogeneous colouring of $\cC(1)$ that preserves the colours of $\cC(0)$.
\item In the second stage, Player A chooses $\cC(2) \varsupsetneq \cC(1)$ and $\cC(2) \subset \cD_j$.
Player B determines an $(\eta,d)$-homogeneous colouring of $\cC(2)$ that preserves the colours of $\cC(1)$.
\item At stage $n$, Player A chooses $\cC(n) \varsupsetneq \cC(n-1)$ and $\cC(n) \subset \cD_j$.
Player B determines an $(\eta,d)$-homogeneous colouring of $\cC(n)$  preserving the colours of $\cC(n-1)$.
\item The game stops at stage $n$ if $\cC(n-1) = \cD_j$, and then Player B is the winner, 
or if there does not exist an $(\eta,d)$-homogeneous colouring of $\cC(n)$ that preserves the colours of $\cC(n-1)$.
In the second case, Player A is the winner.
\end{enumerate}
\paragraph{Defining homogeneous colourings.}
For a collection
$\cc \subset \cd_j $, consider its partition into $d$ subcollections
 $\cc = \cc_1 \cup \ldots \cup \cc_d$.
 Such partition we call colouring of $\cC$.
Colour of  $\Gamma \in \cc$ means the (unique) index $i$ such that
$\Gamma \in \cc_i$.

Let $L \in \cd $, $|L| \geq {1 \over 2^j}$. Denote
\begin{equation}\numer{counting}
\rho(\cc, L) =  |\{ \Gamma \in \cc:  \Gamma \subset L\}|,
\quad\quad
\rho_{i}(\cc, L)  =  |\{ \Gamma \in \cc_i:  \Gamma \subset L\}|.
\end{equation}
\begin{defi}
\numer{c.d1}
Let $\cc \subset \cd_j $, and fix $d \in \bN$, $0<\eta \leq {1 \over 2}$. 
Let   $\cc = \cc_1 \cup \ldots \cup \cc_d$ be some decomposition of $\cc$.
This decomposition is called $(\eta,d)$-homogeneous  colouring of $\cC$ if for each $L \in \cd$, $|L| \geq {1 \over 2^j}$
one of the following holds: 
\begin{itemize}
\item[] 
{\bf Either} $\rho(\cc,L) \leq d$, and then 
\begin{equation}
\numer{hom1}
\rho(\cc_i,L) \leq 1 \quad \text{ for each } \quad 1 \leq i \leq d,
\end{equation}
\item[]
{\bf Or else} $\rho(\cc,L) >d$, and then
\begin{equation}
\numer{hom2}
\eta \max_{1 \leq i \leq d} \rho_i(\cc,L) \leq  \min_{1 \leq i \leq d} \rho_i(\cc,L).
\end{equation}
\end{itemize}
\end{defi}
\paragraph{Remark:} We remark that 
for each (uncoloured) $\cc \subset \cd_j$, $d\in \bN$ and $\eta = {1 \over 2}$ there is always 
a $(\eta,d)$-homogeneous decomposition that can be
obtained as follows: Enumerate the intervals in $\cc$ from left to right, 
and simply put 
$$
\cc_r = \{\Gamma_l \in \cc: l = r \hbox{ mod } d\},\quad\quad 1\le r \leq d.
$$

The problem we treat in this section consists of finding $(\eta,d)$-homogeneous  decompositions
of $\cC \cup \cU $ under the condition that a previously given $(\eta,d)$-homogeneous  decomposition
of $\cC$ is preserved. More precisely,
given an $(\eta,d)$-homogeneous  decompositions of $\cC $ and given $\cU$ disjoint from
$\cC$ we seek to determine a splitting of $\cU$ as 
$$ \cU = \cU _1 \cup \dots \cup \cU_d, $$
which induces an $(\eta,d)$-homogeneous  decomposition of $\cU \cup \cC $ as 
$$ \cU \cup\cC = \cU _1 \cup\cC_1 \cup \dots \cup \cU_d \cup\cC_d .$$
We refer to this question as the decompositon/allocation problem for $\cU . $
Note when $\cU $ contains one element only our problem is just concerned with
allocation. It should be pointed out that the smaller $\cU$ is, the harder it is to find 
a suitable splitting.
Clearly the following possibilities may arise.
\begin{enumerate}
\item The solution to the decomposition/allocation  problem is unique. That is, there exists
just one decomposition of $\cU  $ so that 
$$ \cU _1 \cup\cC_1 , \dots , \cU_d \cup\cC_d $$
is an $(\eta,d)$-homogeneous  decomposition of $\cU \cup \cC .$ 
\item The decompositon/allocation  problem for $\cU  $ does not have a solution.
\item There are conditions on $\cU $ and $\cC$ implying
that the  decompositon/allocation  problem for $\cU  $ has a solution.
Here it is important that any such  condition refers just to 
$\cC$ and not
to the given and fixed decomposition $\{\cC_i \}.$
\end{enumerate}
In this section we address  these three possibilities and examine the transition from one 
case to the next.
\begin{enumerate}
\item We isolate a condition on $\cU$ and $\cC$ (previsibility; see Definition~\ref{c.d2})
implying
that the  decompositon problem for $\cU  $ has a solution. 
See Theorem~\ref{c.t1} which gives rise to winning strategies for Player B.
\item We give  examples where the decomposition/allocation problem for $\cU$ has just one solution. 
Moreover, we give examples (of $\cC$, its decomposition $\{\cC_i\}$ and $\cU$) for which
 the decompostion problem does not have a solution. See Proposition~\ref{e.s1}. This translates to 
 an initial configuration of the game, where Player A has a winning strategy. 
\end{enumerate}

\subsection{Winning Strategies for Player B.}
In the following definition  we isolate a criterion under which Player B can always make his/her move. 
Recall that for a dyadic interval $L \in \cD$, we say that the intervals $L', L'' \in \cD$ are its dyadic successors if $L=L' \cup L''$,
 $|L'|=|L''|={1 \over 2}|L|$.  
\begin{defi}
\numer{c.d2}
Let $\cc, \cu \subset \cd_j$, $\cc \cap \cu = \emptyset$. Let $d \in \bN$.
The collection $\cu$ is called $d$-previsible with respect to the collection $\cc$
 if for every $L \in \cd$ with $|L| \geq {1 \over 2^{j-1}}$ and its dyadic succesors $L', L''$, 
 the following holds: 
  $$
 \rho(\cu\cup \cc,L') < d \quad \hbox{and} \quad  \rho(\cu\cup \cc,L'') \geq d \quad 
 \hbox{implies} \quad  \rho(\cu, L'')=0.
 $$
\end{defi}
Now, we have the following Theorem that imposes
 restrictions on the game -- specifically on the feasible choices of moves for Player A -- under which
Player B can build a winning strategy.

\begin{theor}
\numer{c.t1}
Fix $d \in \bN$ and $\eta$, $0<\eta \leq {1 \over 2}$.
Let $\cc \subset \cd_j$, and let $\{\cc_i, 1 \leq i \leq d\}$ be a fixed $(\eta,d)$-homogeneous  decomposition of $\cc$.
Let $\cu \subset \cd_j$ be $d$-previsible with respect to $\cc$. Then there is a decomposition
$\{\cu_i, 1 \leq i \leq d\}$ of $\cu$ such that $\{\cc_i \cup \cu_i, 1 \leq i \leq d\}$ is an  
$(\eta,d)$-homogeneous  decomposition of $\cc\cup\cu$.
\end{theor}

\proof
 Denote $\ch = \cc \cup \cu$. We are going to define partition of $\cu$ by an inductive argument.
Let $\alpha$ be such that $2^{\alpha }  \leq d < 2^{\alpha+1}$. 
Let us observe that if 
${1 \over 2^j} \leq |L| \leq {1 \over 2^{j-\alpha}}$, then $\rho(\cH, L) \leq 2^\alpha \leq d$.
Thus, if the homogeneity
conditions \eqref{hom1} respectively \eqref{hom2}  are satisfied for $L \in \cd$
with $|L| \geq {1 \over 2^{j-\alpha}}$, then they are satisfied for each $L \in \cd$ with $|L| \geq {1 \over 2^j}$.
Therefore, in our procedure of colouring $\cu$  we consider only $L \in \cd_k$ with $k \leq j-\alpha$.

\smallskip

{\bf I.} $L \in \cd_{j-\alpha}$. Then either  $\rho(\ch,L)<d$ or  $\rho(\ch,L)=d$.

\smallskip

{\bf I.1.} If $\rho(\ch,L)<d$, then intervals from $\cu$ included
in $L$ are left uncoloured.

\smallskip

{\bf I.2.} If $\rho(\ch,L)=d$, then also $\rho(\cc,L) \leq d$, which implies that
$\rho_i(\cc,L) \leq 1$ for each $1 \leq i \leq d$. In such case it is possible to colour
intervals from $\cu$ included in $L$ so that $\rho_i(\ch, L) = 1$
for each $i$, $1 \leq i \leq d$.

\medskip

{\bf II.}   
$L \in \cd_\nu$ with $\nu<j-\alpha$. Then $L = L' \cup L''$
with $L' , L'' \in \cd_{\nu+1}$, $\nu+1 \leq j - \alpha$.
Each interval from $\cc$ or $\cu$ included in $L$ is included in $L'$ or
$L''$, so we have
\begin{eqnarray*}
 \rho(\cc, L) &  = & \rho(\cc, L') + \rho(\cc, L''),
\\
 \rho(\cu, L) &  = & \rho(\cu, L') + \rho(\cu, L''),
\\
 \rho(\ch, L) &  = & \rho(\ch, L') + \rho(\ch, L'').
\end{eqnarray*}

\medskip

\noindent
{\bf Induction hypothesis:} Let $K \in \cd_{\nu+1}.$ 
 If $\rho(\ch,K) <
d$, then intervals from $\cu$ included in $K$ are still uncoloured. 
If $\rho(\ch, K) \geq d$, then all intervals from $\cu$ included in $K$
are coloured, hence the counting parameters $\rho_i(\ch,K)$ are well defined by
\eqref{counting}. They satisfy  $\rho_i(\ch,K) \geq 1$ and
\begin{equation}
\numer{(b3)}
\eta \max_{1 \leq i \leq d} \rho_i(\ch,K) \leq \min_{1 \leq i \leq d} \rho_i(\ch,K).
\end{equation}
\medskip

Then we have two main cases:

\medskip

\noindent
{\bf II.1.} $\rho(\ch,L) < d$. Clearly,  then also $ \rho(\ch, L'),
\rho(\ch, L'') <d$, and intervals from $\cu$ included in $L', L''$
are  uncoloured. 

If $\nu >0$, then leave intervals from $\cu$
included in $L$ still uncoloured.

If $\nu=0$, then $L=[0,1]$, and the induction ends. This means that
$|\ch| < d$, and it is enough to assign elements of $\cu$ to
colours different from colours of elements of  $\cc$.

\medskip

\noindent
{\bf II.2.} $\rho(\ch, L) \geq d$. Then we have next four subcases:

\smallskip

\noindent
{\bf II.2.1.} $\rho(\ch, L') \geq d$ and $\rho(\ch, L'') \geq
d$.
 Then by induction hypotesis all intervals from $\cu$ included in
 $L'$ and $L''$ are already coloured, i.e. all intervals from $\cu$
 included in $L$ are coloured. Moreover, for each $1 \leq i,k \leq d $
$$
\eta \rho_{i}(\ch,L) = \eta  \rho_{i}(\ch,L') + \eta \rho_{i}(\ch,L'')
\leq  \rho_{k}(\ch, L') +  \rho_{k}(\ch, L'') =   \rho_{k}(\ch,
L).
$$
Of course, we have also $\rho_{i}(\ch,L) \geq 1$.

\medskip

\noindent
{\bf II.2.2.} $\rho(\ch, L') < d$ and $\rho(\ch, L'') <
d$. Then by induction hypotesis all intervals from $\cu$ included in
 $L'$ and $L''$ are uncoloured, but the intervals from $\cc$ carry
 their colours. 

Now, we need to colour all  intervals from $\cu$ included in $L', L''$. To simplify notation, let
$m=\rho(\cc, L')$, $n = \rho(\cc, L'')$, $x=\rho(\cu,L')$,
$y=\rho(\cu,L'')$. We have $0 \leq m, n \leq d-1$, $0 \leq m+x,n+y
\leq d-1$ and $d \leq \rho(\ch, L) = m+x+n+y \leq 2(d-1)$.

First consider the case $m+n <d$. Then
 $0 \leq \rho_{i}(\cc,L) \leq
1$ for each $i$. 
For simplicity, assume that intervals from $\cc$ included in $L'$
have colours $1, \ldots, m$, and intervals from $\cc$ included in
$L''$ have colours $m+1, \ldots , m+n$.
Now, we colour intervals from $\cu$.
First, colour intervals from $\cu$ included in $L'$ using 
colours $m+n+1, \ldots, d$, and then, if necessary (i.e. $x > d - (m+n)$), continuing with $x
- (d - (m+n))$
colours from  $m+1, \ldots, m+n$; since $m+x <d$, in this way we assign
colours to all intervals from $\cu$ included in $L'$.
 Next, we assign colours to
intervals from $\cu$ included in $L''$. If $m+n+x<d$, then assign
first colours $m+n+x+1, \ldots, d$, then continue with colours $1,
\ldots, m$, and then if necessary with colours $m+n+1, \ldots, m+n+x$. 
If $m+n+x \geq d$, then just choose $y$ different colours from $1,
\ldots,m$ and $m+n+1, \ldots, d$. With such colouring of intervals
from $\cu$ included in $L'$ and $L''$ we find that 
both
$\rho_{i}(\ch,L') \leq 1$ and
  $\rho_{i}(\ch,L'') \leq 1$.
 This implies that for each $K \subset L'$ or $K\subset
L''$ we have $\rho(\ch,K) < d$ and $\rho_{i}(\ch, K) \leq 1$.
Moreover, we get $1 \leq \rho_{i}(\ch,L) \leq 2$, which implies
$$
\eta \max_{1 \leq i \leq d} \rho_i(\ch,L) \leq {1 \over 2} \max_{1 \leq i \leq d} \rho_i(\ch,L)
\leq \min_{1 \leq i \leq d} \rho_i(\ch,L).
$$

It remains to consider the case $m+n \geq d$. 
Then the homogeneity assumption on the decomposition of $\cc$ -- \eqref{hom1} for $L',L''$ and \eqref{hom2} for $L$ --  implies
$1 \leq \rho_i(\cc,L)  = \rho_i(\cc,L') + \rho_i(\cc,L'') \leq 2$.
 For simplicity, assume
that intervals from $\cc$ included in $L'$ have colours $1, \ldots,
m$ and intervals from  $\cc$ included in $L''$ have colours $m+1, \ldots, d$ and
$1, \ldots, m+n-d$. To colour intervals from $\cu$ included in $L'$
choose $x$ colours from $m+1, \ldots, d$. To colour intervals from
$\cu$ included in $L''$ choose $y$ colours from $m+n-d+1, \ldots , m$.
This is possible since $m+x < d$ and $n+y<d$. Observe that in this way
we get $0 \leq \rho_{i}(\ch,L'), \rho_{i}(\ch,L'') \leq 1$ and $1
\leq \rho_{i}(\ch,L) \leq 2$. Therefore, for each $K \subset L'$ or $K \subset L''$
we have $\rho_i(\ch,K) \leq 1$, while for $L$ we have
$$
\eta \max_{1 \leq i \leq d} \rho_i(\ch,L) \leq {1 \over 2} \max_{1 \leq i \leq d} \rho_i(\ch,L)
\leq \min_{1 \leq i \leq d} \rho_i(\ch,L).
$$

\medskip

\noindent
{\bf II.2.3.} $\rho(\ch, L') < d$ and $\rho(\ch, L'') \geq
d$. 
Then by induction hypotesis all intervals from $\cu$ included in
 $L'$  are uncoloured,  but the intervals from $\cc$ included in
 $L'$ carry their colours.
Since $\cu$ is $d$-previsible with respect to $\cc$, we have $\rho(\cu, L'') = 0$.
Therefore, 
$\rho(\ch,L'') = \rho(\cc, L'')$, and by  condition \eqref{hom2} of the 
$(\eta,d)$-homogeneity for $\cc$, we get
$\rho_i(\ch,L'') = \rho_i(\cc,L'') \geq 1$ and $\rho_{i}(\ch,L'')$ satisfy \rf{(b3)}.

\smallskip

If $\rho(\cu, L') = 0$ as well, then all intervals from $\ch$
included in $L$ come from $\cc$, and there is nothing to do.

\smallskip

Let $\rho(\cu,L') =x > 0$. We need to colour $x$ intervals from
$\cu$ included in $L'$. To simplify notation, let $m = \rho(\cc,
L')$. Note that $1 \leq m+x <d$. Let $S = \{i: \rho_{i}(\cc,L')
=1\}$ and $T= \{i:\rho_{i}(\cc,L') =0\}$. Let $t_1, \ldots t_{d-m}$
be an ordering of $T$ such that 
\begin{equation}
\numer{(c4)}
\rho_{t_1}(\cc,L'') \leq \ldots \leq \rho_{t_{d-m}}(\cc,L'').
\end{equation}
Since $x < d-m$, there are more colours in $T$ than intervals in $\cU$ that are included in $L'$.
Now  attach the colours ${t_1}, \ldots, {t_x}$, bijectively, to intervals in $\cU$  contained in $L'$.
Then $\rho_{i}(\ch,L') \leq 1$. 

It remains  to check that
$\rho_{i}(\ch,L)$ satisfy \rf{(b3)}. By assumption on partition of
$\cc$ we have
\begin{equation}
\numer{(c1)}
\eta \max_{1 \leq i \leq d} \rho_{i}(\cc, L'') 
\leq  \min_{1 \leq i \leq d} \rho_{i}(\cc, L''), 
\end{equation}
\begin{equation}
\numer{(c2)}
\eta \max_{1 \leq i \leq d} \rho_{i}(\cc, L) 
\leq  \min_{1 \leq i \leq d} \rho_{i}(\cc, L).
\end{equation}
Moreover,
\begin{equation}
\numer{(c3)}
\begin{cases}
\rho_{i}(\ch,L) = \rho_{i}(\cc, L) = \rho_{i}(\cc, L'') +1 
& \text{ for } i \in  S,
\cr \cr
\rho_{i}(\ch,L) = \rho_{i}(\cc, L) +1 = \rho_{i}(\cc, L'') + 1
& \text{ for } i=t_1, \ldots, t_x,
\cr\cr
\rho_{i}(\ch,L) = \rho_{i}(\cc, L) = \rho_{i}(\cc, L'') 
& \text{ for } i=t_{x+1}, \ldots t_{d-m}.
\end{cases}
\end{equation}
Observe that if  $\max_{i} \rho_{i}(\ch,L) = \rho_{k}(\ch,L)$, then
$k \in S$ or $k \in T$; in case $k \in T$ we have  $k = t_x$ or
$k=t_{d-m}$, because of ordering \rf{(c4)}. 
If  $k \in S$, then
\rf{(b3)} is satisfied for $L$ and $\ch$ because of \rf{(c2)} and the first line of
\rf{(c3)}. If $k=t_{d-m}$, then
\rf{(b3)} is satisfied for $L$ and $\ch$ because of \rf{(c1)} and the last line of
\rf{(c3)}. If $k = t_x$ and $\rho_{t_x}(\ch,L) > \rho_{t_{d-m}}(\ch,L)$   then we check inequality
\begin{equation}
\numer{(c5)}
\eta \rho_{t_x}(\ch,L) \leq  \rho_{i}(\ch,L) \quad {\rm for} \quad 1
\leq i \leq d.
\end{equation}
For $i \in S$  inequality \rf{(c5)} is satisfied because of \rf{(c1)} and the first two
lines of \rf{(c3)}. For $i= t_1, \ldots, t_x$   inequality \rf{(c5)} is satisfied
because of \rf{(c2)} and the second line of \rf{(c3)}. When
$\rho_{t_x}(\ch,L) > \rho_{t_{d-m}}(\ch,L)$, then the two last lines
of \rf{(c3)} and the ordering \rf{(c4)} imply $\rho_{t_x}(\cc,L'') = \rho_{t_{x+1}}(\cc,L'')=
 \ldots = \rho_{t_{d-m}}(\cc,L'')$. This implies that inequality
 \rf{(c5)} is satisfied, even  with  ${1 \over 2}$ on the left-hand-side, for
 $i=t_{x+1}, \ldots t_{d-m}$.

\bigskip

\noindent
{\bf II.2.4.} $\rho(\ch, L') \geq d$ and $\rho(\ch, L'') <
d$. This case is analogous to II.2.3.

\medskip

This completes the proof of Theorem \ref{c.t1}.
\endproof

\subsection{Winning Strategies for Player A.}
Here we  analyze the role of the  previsibility assumption in Theorem~\ref{c.t1}.
We do this by defining an initial configuration of the two-person game so that Player A has a strategy to win in exactly $n$
moves. 
This corresponds to a sequence of examples for which the 
decomposition/allocation problem has a unique  solution and  a related example for which 
the decomposition/allocation problem is without solution.

We start with $\cC(0)$, its initial decomposition and $\cU(0)$ in such a way that the decomposition  problem for
$\cU(0)$ has just one solution. This uniquely determined solution defines the decomposition for 
$\cC(1)= \cU(0) \cup \cC(0).$ This and the given $\cU(1) $ determines a decomposition problem 
for which we will see that it again has only one solution.
This solution in turn determines a splitting of  $\cC(2) = \cU(1) \cup \cC(1)$ which again leads to a 
decomposition problem with a unique solution. This will go on until we reach $\cC(n-1) $ 
and its decomposition that has been determined uniquely by $\cC(0) $ its initial decomposition
and by our choice of $\cU(0),\dots ,\cU(n-2).$ 
Then we change the situation and choose the collection  
$\cU(n-1)$ that forces the decomposition problem in $\cU(n-1)\cup \cC(n-1)$ to be without solution.

Throughout this section we take 
$d = 2^a$, $a \in \bN$, and $\eta = {1 \over n}$ with $n \in \bN$ and $j \ge n+a +1 .$
\begin{prop}\numer{e.s1}
There exist
$$
\cC(0) \subset \cD_j
$$
with $(\eta,d)$-homogeneous decomposition 
$$\cc_1(0), \ldots , \cc_d(0)$$
and an increasing chain of collections
$$ \cc(0) \subset \cc(1) \subset \ldots \subset \cc(n) \subset \cd_j$$
so that for
$$
\cu(k) = \cc(k+1) \setminus \cc(k)\quad \text{ with } \quad 1 \leq k \leq n-1,
$$
the following conditions hold:
\begin{itemize}
\item[(A)] Stage $0$.
There exists exactly one splitting of $\cu(0)$ as
$
\cu_1(0), \ldots, \cu_d(0)
$
so that
$$
\cu_1(0)\cup \cc_1(0), \ldots, \cu_d(0)\cup \cc_d(0)
$$
is an $(\eta,d)$-homogeneous decomposition of $\cu(0) \cup \cc(0)$, hence of  $\cc(1)$.

\item[(B)] Stage $k$, $1 \leq k \leq n-2$,
let 
$$
\cc_1(k), \ldots, \cc_d(k)
$$
be the unique $(\eta,d)$-homogeneous decomposition of $\cc(k)$, obtained at stage $k-1$.
There exists exactly one splitting of $\cu(k)$ as
$
\cu_1(k), \ldots, \cu_d(k)
$
so that
$$
\cu_1(k)\cup \cc_1(k), \ldots, \cu_d(k)\cup \cc_d(k)
$$
is an $(\eta,d)$-homogeneous decomposition of $\cu(k) \cup \cc(k)$, hence of  $\cc(k+1)$.

\item[(C)] Stage $n-1$. Let 
$$
\cc_1(n-1), \ldots, \cc_d(n-1)
$$
be the unique $(\eta,d)$-homogeneous decomposition of $\cc(n-1)$, obtained at stage $n-2$.
There does not exist a splitting of $\cu(n-1)$ as
$
\cu_1(n-1), \ldots, \cu_d(n-1)
$
so that
$$
\cu_1(n-1)\cup \cc_1(n-1), \ldots, \cu_d(n-1)\cup \cc_d(n-1)
$$
is an $(\eta,d)$-homogeneous decomposition of $\cc(n) = \cu(n-1) \cup \cc(n-1)$.
\end{itemize}
\end{prop}
\proof  Observe that
for each $j$, the testing  levels for $\cc \subset \cd_j$ are
$\cd_{j-a}, \cd_{j-a-1}, \ldots, \cd_1, \cd_0$.
Since $j \geq n+a+1$,  there are at least $n+2$ testing levels.
Take a chain of dyadic intervals
 $$L_1 \subset L_2 \subset \ldots \subset L_{n+2},
\quad L_i \in \cd_{j-a-i+1}.$$
Then $|L_i| = {1 \over 2} |L_{i+1}|$, and let $P_i$ be the
dyadic brother of $L_i$ in $L_{i+1}$, $ i = 1, \ldots n+1$. Thus $P_i=L_{i+1} \setminus L_i$.

Now, take two sets of intervals from $\cd_j$:
$$I_1, \ldots, I_{d-1} \in \cd_j \text{ such that } I_i \subset L_1
\text{ for
each } i= 1, \ldots,  d-1,$$
$$J_1, \ldots , J_{n+1} \in \cd_j \text{ such that }
J_i \subset P_i \text{ for each }
i = 1, \ldots , n+1.
$$
Consider the following sequence of collections:
\begin{eqnarray*}
\cc(0) & = & \{ I_1, \ldots, I_{d-1}\} \cup \{ J_{n+1}\},
\\
\cc(k)&  = & \{I_1, \ldots, I_{d-1}\} \cup \{J_{n-k+1}, \ldots, J_{n+1}\}.
\end{eqnarray*}
Observe that
$$
\cu(k)  =  \{J_{n-k}\}, \quad k = 0, \ldots, n-1,
$$
and
\begin{equation}
\numer{j08.e1}
J_i \in \cc(n+1-i),
\end{equation}
and since our chain of collections $\cC(0), \ldots, \cC(n)$ is increasing, $J_i$ is also contained in $\cc(k)$ with
$k\geq n+1-i$.

\paragraph{Initialization -- verification of (A).}  Consider possible colourings of $\cc(0)$.
Take $L_{n+2}$ as a testing interval. Observe that $\rho(\cc(0))
= \rho(\cc(0),L_{n+2})=d$, so  if we want to have $(\eta,d)$-homogeneity, 
we must have \eqref{hom1} and therefore
$ \rho_{i}(\cc(0), L_{n+2}) =1 $ for each $i=1, \ldots,d$.
 Without loss of generality we can assume
that $J_{n+1}$ has colour 1, and each $I_i$ has colour $i+1$, $i=1, \ldots, d-1$.
Therefore for $\cc(0)$ and each testing interval $L \subset L_{n+2}$ we have
$\rho_i(\cc(0),L) \leq 1$ for each $1 \leq i \leq d$.

\paragraph{The basic observation.} Our example is based on iterating systematically the following basic observation.
Let $k \leq n$.
Assume that $\cc(k)$ has an $(\eta,d)$-homogeneous decomposition as
$$
\cc_1(k), \ldots , \cc_d(k),
$$
so that
$$
\cc_1(0)\subset \cc_1(k), \ldots , \cc_d(0) \subset \cc_d(k).
$$
Then necessarily 
\begin{equation}
\numer{jan08.e2}
J_{n-k+1} \text{ must have colour } 1.
\end{equation}
\paragraph{Verification of \eqref{jan08.e2}.}
We know already that  $J_{n+1}$ has to have colour $1$. To check the claim for $J_{n-k+1}$, $k=1, \ldots, n$
we consider the pair of collections $\cc(0) \subset \cc(k)$:
$$
 \cc(k)  =   \cc(0) \cup \{J_{n-k+1}, \ldots, , J_{n}\}.
$$
and testing interval $L_{n-k+2}$. Elements of $\cc(0)$ included in $L_{n-k+2}$ are
$I_1, \ldots, I_{d-1}$. In addition, $J_{n-k+1} \subset P_{n-k+1} \subset L_{n-k+2}$, 
while
$J_{n-k+2}, \ldots, J_n \not \subset L_{n-k+2}$.
Therefore we have
$$
\rho(\cc(0), L_{n-k+2}) = d-1, \quad
\rho(\cc(k), L_{n-k+2}) = d,$$
$$\rho_{1}(\cc(0), L_{n-k+2}) = 0 \quad
{\rm and} \quad  \rho_{i}(\cc(0), L_{n-k+2}) = 1 \quad {\rm  for} \quad  i=2, \ldots, d.$$
Therefore, \eqref{hom1} of the $(\eta,d)$-homogeneity condition for 
$\cc(k)$ implies that $J_{n-k+1}$ is of colour 1.

\paragraph{Verificaton of (B).}
Recall that $0 \leq k \leq n$
$$
\cc(k) =  \{I_1, \ldots, I_{d-1}\} \cup \{J_{n-k+1}, \ldots J_{n+1}\}.
 $$
Moreover, by \eqref{jan08.e2},  the only possible $({1 \over n},d)$-homogeneous decomposition of 
$\cc(k)$ is
$$
\cc_1(k) = \{J_{n-k+1}, \ldots J_{n+1}\}, \quad \cc_i(k) = \{ I_{i-1} \} \quad \hbox{for} \quad 2 \leq i \leq d
$$
Let's check that for $0 \leq k \leq n-1$, the above decomposition of $\cc(k)$ is indeed $({1 \over n},d)$-homogeneous.
We present the detailed proof for $k=n-1$, since  the cases $ k \leq n-1$ are fully analogous.

First, take as a testing  interval $L_{s}$, $s=3, \ldots n+2$.
Then  elements of $\cc(n-1)$  included in $L_{s}$
are $I_1, \ldots, I_{d-1}$ and $J_2, \ldots, J_{s-1}$.
Therefore $$\rho(\cc(n-1), L_s) = s+d-3,$$ and
$$
\rho_{1}(\cc(n-1), L_s) = s-2,
\quad \rho_{i}(\cc(n-1), L_s) = 1 \quad {\rm for} \quad i=2, \ldots d.
$$
Therefore 
$${1 \over n} \max_{1 \leq i \leq d} \rho_i(\cc(n-1), L_s)
\leq \min_{1 \leq i \leq d} \rho_i(\cc(n-1), L_s), \quad s = 3, \ldots, n+2.$$

Next take as a testing interval $L_2$.
Then  elements of $\cc(n-1)$  included in $L_{2}$
are $I_1, \ldots, I_{d-1}$, so  $\rho(\cc(n-1), L_2) = d-1$,
 $$
\rho_{1}(\cc(n-1), L_2) = 0,
\quad \rho_{i}(\cc(n-1), L_2) = 1 \quad {\rm for} \quad i=2, \ldots d.
$$
Therefore $L_2$ also satisfies \eqref{hom1}  of the $({1 \over n},d)$-homogeneity condition for $\cc(n-1)$. 
Consequently, $L_1, P_1 \subset L_2$ also satisfy these conditions.

Finally, take as a testing interval $P_k$, $k=2,\ldots, n+1$.
The only element of $\cc(n-1)$ included in $P_k$ is $J_k$,  so $\rho(\cc(n-1),P_k)=1$, and more precisely
 $$
\rho_{1}(\cc(n-1), P_k) = 1,
\quad \rho_{i}(\cc(n-1), P_k) = 0 \quad {\rm for} \quad i=2, \ldots d.
$$
Thus, $P_k$ (and consequently, each testing interval included in $P_k$) satisfies \eqref{hom1} 
of the $({1 \over n},d)$-homogeneity condition for $\cc(n-1)$.

\paragraph{Verification of (C).}
Consider $\cc(n-1)$ and $\cc(n) = \cc(n-1) \cup \cu(n-1)$.
Recall that
$$
 \cc(n) =  \{I_1, \ldots, I_{d-1}\} \cup \{J_1, J_2, \ldots J_{n+1}\}.
$$
Take    $L_{n+2}$ as a testing interval. All intervals from $\cc(n)$ are included in $L_{n+2}$, and the colouring yields
 $$
\rho_{1}(\cc(n), L_{n+2}) = n+1,
\quad \rho_{i}(\cc(n), L_{n+2}) = 1 \quad {\rm for} \quad i=2, \ldots d.
$$
For $\cc(n)$ and $L_{n+2}$ we have to consider \eqref{hom2}   of the $({1 \over n},d)$-homogeneity condition.
But the above formulae  mean that for $\cc(n)$ and testing interval $L_{n+2}$, the condition \eqref{hom2}  
 is satisfied
with $\eta' = {1 \over n+1}$, but not with $\eta = {1 \over n}$.\endproof

\paragraph{Remark.} For $0 \leq k \leq n-1$, the collection $\cU(k)$ is not previsible with respect to $\cC(k)$.
Nevertheless, the colouring problem has a solution for $0 \leq k \leq n-2$.

\nocite{gmp} 
\nocite{ms}
\nocite{haydon-2007}

\bibliographystyle{abbrv}

\bibliography{isojan09sept091}

\begin{tabular}{lll}
Anna Kamont & \hspace{4.5cm} &Paul F.X. M\"uller
\\
Institute of Mathematics & & Department of Analysis
\\
Polish Academy of Sciences & &J.Kepler University
\\
ul. Abrahama 18, 81-825 Sopot & &A-4040 Linz
\\
Poland & & Austria
\\
A.Kamont@impan.gda.pl & & pfxm@bayou.uni-linz.ac.at
\end{tabular}
\end{document}